\documentclass[a4paper,12pt]{article}
\usepackage{amsthm}
\usepackage{amsmath}
\usepackage{amssymb}
\usepackage{latexsym}
\usepackage{url}
\usepackage{color}

\usepackage{tikz}
\usetikzlibrary{intersections, calc, arrows.meta}

\usepackage{comment}
\numberwithin{equation}{section}

\def\R{{\bf R}}

\def\N{{\bf N}}

\def\d{\displaystyle}
\def\e{{\varepsilon}}

\def\wt{\widetilde}

\def\p{\partial}
\def\v#1{\mbox{\boldmath $#1$}}

\newcommand{\al}{\alpha}
\newcommand{\be}{\beta}

\newcommand{\x}{\xi}

\newcommand{\EQS}[1]{\begin{align} #1 \end{align}}

\newcommand{\LR}[1]{{\langle {#1} \rangle }}

\newtheorem{thm}{Theorem}[section]

\newtheorem{lem}{Lemma}[section]

\newtheorem{rem}{Remark}[section]

\title{The lifespan estimates of classical solutions of one dimensional semilinear wave equations\\ with characteristic weights}
\author{
Shunsuke Kitamura,
\footnote{Doctor course, Mathematical Institute, Tohoku University, Aoba, Sendai 980-8578, Japan. e-mail: shunsuke.kitamura.s8@dc.tohoku.ac.jp.}
\  \
Hiroyuki Takamura,
\footnote{Mathematical Institute, Tohoku University, Aoba, Sendai 980-8578, Japan. e-mail: hiroyuki.takamura.a1@tohoku.ac.jp.}\\
\  \
Kyouhei Wakasa
\footnote{Department of Creative Engineering, National Institute of Technology, Kushiro College, 2-32-1 Otanoshike-Nishi, Kushiro-Shi, Hokkaido 084-0916, Japan. e-mail: wakasa@kushiro-ct.ac.jp.
}
}
\date{
\[
\begin{array}{llll}
\mbox{\footnotesize{\bf Keywords:}}
& \mbox{\footnotesize semilinear wave equation, one space dimension, classical solution,}\\
& \mbox{\footnotesize lifespan}\\
\mbox{\footnotesize{\bf MSC2020:}}
& \mbox{\footnotesize Primary 35L71, Secondary 35B44}\\
\end{array}
\]
}
\begin{document}
\maketitle
\begin{abstract}
In this paper, we study the lifespan estimates of classical solutions for semilinear wave equations
with characteristic weights and compactly supported data in one space dimension.
The results include those for weights by time-variable, but exclude those for weights by space-variable in some cases.
We have interactions of two characteristic directions.
\end{abstract}

%%%%%%%%%%%%%%%%%%%%%%%%%%%%%%%%%%%%%%%%%%%%%%%%%%
%%%%%%%%%%%%%%%%%%%%% SECTION1 %%%%%%%%%%%%%%%%%%%%%%%
%%%%%%%%%%%%%%%%%%%%%%%%%%%%%%%%%%%%%%%%%%%%%%%%%%
%\begin{tikzpicture}
% \draw[->,>=stealth,semithick] (-5,0)--(5,0)node[above]{$a$}; %x軸
% \draw[->,>=stealth,semithick] (0,-5)--(0,5)node[right]{$b$}; %y軸
% \draw (0,0)node[below left]{O}; %原点
% \draw (-1,0)node[below]{$-1$}; %点(-1,0)
% \draw (0,1)node[right]{$1$}; %点(0,1)
% \draw[very thick,domain=0:4.5] plot(\x,-\x)node[right]{$b=-a$};
%\end{tikzpicture}

\section{Introduction}
In this paper, we focus on the study of a model equation for the purpose of extending the general theory of nonlinear wave equations. Before stating our target, we first
overview the general theory for nonlinear wave equations in one space dimension. Consider the following Cauchy problem:
\EQS{\label{GIVP}
\left\{
\begin{array}{ll}
\d  u_{tt}- u_{xx}=H(u,u_t,u_x,u_{xx},u_{xt}) &\quad \mbox{in}\ \R\times (0,\infty),\\
 u(x,0)= \e f(x) , \quad u_t(x,0)= \e g(x) &\quad \mbox{for}\ x\in\R,
\end{array}
\right.
}
where $f,g\in C_0^{\infty}(\R)$ and $\e>0$ is a small parameter. 
Let $\wt{\lambda}=(\lambda;\ (\lambda_i), i=0,1; \ (\lambda_{ij}), i,j,=0,1, i+j\ge1).$
Assume that $H=H(\wt{\lambda})$ is a sufficiently smooth function with 
\[
H(\wt{\lambda})=O(|\wt{\lambda}|^{1+\alpha})
\]
in a neighborhood of $\wt{\lambda}=0$, where $\alpha \in \N$. Let us define the lifespan $\wt{T}(\e)$ as the maximal existence time of the classical solution of (\ref{GIVP}) with arbitrary fixed non-zero data.
The general theory for the problem (\ref{GIVP}) is to express the lower bound of the lifespan according to $\alpha$ and the initial data. 
Li, Yu and Zhou \cite{LYZ} have constructed the general theory for this problem: 
\[
\wt{T}(\e)\ge 
\left\{
\begin{array}{lll}
c\e^{-\alpha/2} & \mbox{in general case,}\\
c\e^{-\alpha(1+\alpha)/(2+\alpha)} & \mbox{if}\ \d \int_{\R}g(x)dx = 0,\\
c\e^{-\alpha} & \mbox{if}\ \p_{u}^{\beta}H(0)=0\ \mbox{for}\ 1+\alpha\le \forall \beta \le 2\alpha.
\end{array}
\right.
\]
It should be noted that Morisawa, Sasaki and Takamura \cite{MST} 
point out that there is a possibility to improve the above general theory by studying $H=|u_t|^p+|u|^q$ so-called \lq\lq combined effect'' nonlinearity. 
For the detail, see the end of Section 2 in \cite{MST}. 

We set $H=|u|^p$ in (\ref{GIVP}) which is a model to ensure the optimality of the general theory. Note that Kato \cite{Kato80} showed the blow-up result for any $p>1$.
Zhou \cite{Z92} has obtained the following lifespan estimates for $p>1$:
\begin{equation}
\label{Zhou-est}
\wt{T}(\e)\sim
\left\{
\begin{array}{lll}
\e^{-(p-1)/2} &\mbox{if}&\d \int_{\R} g(x)dx \neq0,\\
\e^{-p(p-1)/(p+1)} &\mbox{if}&\d \int_{\R} g(x)dx =0.
\end{array}
\right.
\end{equation}
Here, the definition of $\wt{T}(\e)\sim A(C,\e)$ is given as follows: there exist positive constants $C_1$ and $C_2$ 
which are independent of $\e$ satisfying $A(C_1,\e)\le \wt{T}(\e)\le A(C_2,\e)$.
The differences between the lifespan estimates of (\ref{Zhou-est}) come from Huygens' principle which holds
if the total integral of the initial speed $g$ vanishes. See Lemma \ref{lm:huygens} below.

Our motivation is to extend the general theory to the non-autonomous non-linear terms: 
\[
H = H(x,t,u,u_t,u_x,u_{xx},u_{xt}).
\]
In order to look for the assumptions on $x$ and $t$ in $H$, we shall investigate the following model equations in this paper:
\EQS{\label{IVP}
\left\{
\begin{array}{ll}
\d  u_{tt}- u_{xx}=F(x,t)|u|^p &\quad \mbox{in}\ \R\times (0,\infty),\\
 u(x,0)= \e f(x) , \quad u_t(x,0)= \e g(x) &\quad \mbox{for}\ x\in\R,
\end{array}
\right.
}
where $p>1$ and $F\in C^1(\R\times (0,\infty))$. Our purpose is to get
the lifespan estimates for the problem (\ref{IVP}) with the characteristic weights:
\begin{equation}
\label{characteristic}
F(x,t)=\LR{t+\LR{x}}^{-(1+a)}\LR{t-\LR{x}}^{-(1+b)},
\end{equation}
where $\LR{x}:=(1+x^2)^{1/2}$ for $x\in \R$ and $a,b\in\R$. 
The meaning of the \lq\lq characteristic'' is originally used as $t+|x|$ and $t-|x|$. We have changed it to (\ref{characteristic}) by using 
$\LR{\cdot}$, because we need the regularity to get the classical solutions 
and to avoid the singularity when the power of the exponent is negative. However, this modification is not essential because of Lemma \ref{weight_equivalent} below.  
Let us denote by $T(\e)$ the maximal existence time of the classical solution of (\ref{IVP}) with arbitrary fixed non-zero data. 
Then, we have the following main results:
\begin{equation}
\label{main_global}
T(\e)=\infty\quad \mbox{for}\ a+b>0\ \mbox{and} \ a>0,
\end{equation}
\begin{equation}
\label{main_lifespan1}
T(\e)\sim
\left\{
\begin{array}{lllll}
\exp(\e^{-(p-1)}) &\mbox{for}&  a+b=0\ \mbox{and}\ a>0, \ \mbox{or}\ a=0\ \mbox{and} \ b>0, \\ 
\exp(\e^{-(p-1)/2}) &\mbox{for}& a=b=0, \\
\e^{-(p-1)/(-a)} &\mbox{for}& a<0 \ \mbox{and}\ b>0,\\
\phi_1^{-1}(\e^{-(p-1)})&\mbox{for}& a<0\ \mbox{and}\ b=0,\\
\e^{-(p-1)/(-a-b)}&\mbox{for}& a+b<0\ \mbox{and}\ b<0\\
\end{array}
\right.
\end{equation}
if 
\[
\int_{\R}g(x)dx\neq0,
\]
where $\phi_1^{-1}$ is an inverse function of $\phi_1$ defined by
\begin{equation}
\label{phi1}
\phi_1(s)=s^{-a}\log(2+s),
\end{equation}
and
\begin{equation}
\label{main_lifespan2}
T(\e)\sim
\left\{
\begin{array}{lllll}
\exp(\e^{-(p-1)}) &\mbox{for}&  a=0\ \mbox{and}\ b>0,\\
\exp(\e^{-p(p-1)}) &\mbox{for}& a+b=0\ \mbox{and} \ a>0, \\  
\exp(\e^{-p(p-1)/(p+1)}) &\mbox{for}& a=b=0, \\
\e^{-(p-1)/(-a)} &\mbox{for}& a<0 \ \mbox{and}\ b>0,\\
\psi_1^{-1}(\e^{-p(p-1)})&\mbox{for}& a<0\ \mbox{and}\ b=0,\\
\e^{-p(p-1)/(-pa-b)} &\mbox{for}& a<0\ \mbox{and}\ b<0,\\
\psi_2^{-1}(\e^{-p(p-1)})&\mbox{for}& a=0\ \mbox{and}\ b<0,\\
\e^{-p(p-1)/(-a-b)} &\mbox{for}& a+b<0 \ \mbox{and}\ a>0
\end{array}
\right.
\end{equation}
if
\[
\int_{\R}g(x)dx=0,
\]
where $\psi_1^{-1}$ and $\psi_2^{-1}$ are inverse functions of $\psi_1$ and $\psi_2$ respectively defined by 
\begin{equation}
\label{psi1_psi2}
\psi_1(s)=s^{-pa}\log(2+s)\ \mbox{and} \ \psi_2(s)=s^{-b}\log^{p-1}(2+s).
\end{equation}

We explain the background to considering the form of the weight function $F$ in (\ref{characteristic}). Firstly, the pointwise estimates of the wave equations have 
a characteristic factor such as (\ref{characteristic}). A natural question arises how the lifespan estimates change as compared with those for non-weighted case, (\ref{Zhou-est}). 
Secondly, our equations include some damped wave equations which were treated by many previous works. 
For example, let us consider the following Cauchy problem for the nonlinear damped wave equations: 
\EQS{\label{scaledampedIVP}
\left\{
\begin{array}{ll}
\d v_{tt}- v_{xx}+\frac{2}{1+t} v_t=|v|^p &\quad \mbox{in}\ \R\times (0,\infty),\\
 v(x,0)= \e f(x) , \quad v_t(x,0)= \e g(x) &\quad \mbox{for}\ x\in\R.
\end{array}
\right.
}
Then, if we set $u(x,t)=(1+t)v(x,t)$ (Liouville transform), we have
\EQS{\label{timeweightIVP}
\left\{
\begin{array}{ll}
\d u_{tt}- u_{xx}=\frac{|u|^p}{(1+t)^{p-1}} &\quad \mbox{in}\ \R\times (0,\infty),\\
 u(x,0)= \e f(x) , \quad  u_t(x,0)= \e\{f(x)+g(x)\} &\quad \mbox{for}\ x\in\R.
\end{array}
\right.
}
When $p>3$, D'Abbicco \cite{DABI} showed that the energy solution of (\ref{scaledampedIVP}) exists globally in time.
On the other hand, Wakasugi \cite{Wakasugi14} has obtained the blow-up results of the problem (\ref{scaledampedIVP}) for $1<p\le3$. 
We note that they treated more general damping terms $\mu v_t/(1+t) $ $(\mu>0)$ in (\ref{scaledampedIVP}).
For the lifespan estimates in (\ref{timeweightIVP}), Wakasa \cite{wak16} obtained
\begin{equation}
\label{wak-est}
\begin{array}{c}
T(\e)\sim
\left\{
\begin{array}{lll}
\e^{-(p-1)/(3-p)} &\mbox{for}&\d 1<p<3,\\
\exp(\e^{-(p-1)}) &\mbox{for}&\d p=3 \\
\end{array}
\right. 
\\
\mbox{if}\ \d \int_{\R}\{f(x)+g(x)\}dx \neq0
\end{array}
\end{equation}
and Kato, Takamura and Wakasa \cite{KTW} obtained
\begin{equation}
\label{KTW-est}
\begin{array}{c}
T(\e)\sim
\left\{
\begin{array}{lll}
\e^{-p(p-1)/(1+2p-p^2)} &\mbox{for}&\d 1<p<2,\\
b(\e) &\mbox{for}&\d p=2,\\
\e^{-p(p-1)/(3-p)} &\mbox{for}&\d 2<p<3,\\
\exp(\e^{-p(p-1)}) &\mbox{for}&\d p=3
\end{array}\right. \\
\mbox{if}\ \d \int_{\R}\{f(x)+g(x)\}dx =0,\end{array}\end{equation}
where $b=b(\e)$ is a positive number satisfying $\e^2 b \log(1+b)=1$. Here, we note that the critical exponent 3  which is a threshold between the global-in-time existence and the blow-up in finite time is the so-called Fujita exponent in one space dimension. 
%Moreover, the lifespan estimate in the case of $1<p<2$ in (\ref{KTW-est}) is same as that of the one of the 3-dimensional semillinear wave equations with the subcritical nonlinearity. 
Remarkably, our lifespan estimates (\ref{main_lifespan2}) with $a=p-2$ and $b=-1$ coinside with (\ref{wak-est}) and (\ref{KTW-est}). 
This is because, $1+t$ is equivalent to $\LR{t+\LR{x}}$ 
by the finite propagation speeds. See (\ref{support_sol}) and Lemma \ref{weight_equivalent} below.

%\EQS{\label{IVP}
%\left\{
%\begin{array}{ll}
%\d  u_{tt}- u_{xx}=\frac{|u|^p}{\LR{t+\LR{x}}^{1+a}\LR{t-\LR{x}}^{1+b}} &\quad \mbox{in}\ \R\times (0,\infty),\\
% u(x,0)= \e f(x) , \quad u_t(x,0)= \e g(x) &\quad \mbox{for}\ x\in\R,
%\end{array}
%\right.
%}
%where $\LR{x}=(1+x^2)^{1/2}$ for $x\in \R$, $a,b\in\R$ and $p>1$. Here, the initial data satisfies  and $\e>0$ is a small parameter. 

We next mention the results to the case where $F$ in (\ref{IVP}) is a spatial weight only as $F=\LR{x}^{-(1+b)}$ with $b\in \R$. 
Kitamura, Morisawa and Takamura \cite{KMT} have obtained the lifespan estimates,
\begin{equation}
\label{lifespan_non-zeroKMT}
T(\e)\sim
\left\{
\begin{array}{ll}
\e^{-(p-1)/(1-b)} & \mbox{for}\ b<0,\\
\phi^{-1}(\e^{-(p-1)}) & \mbox{for}\ b=0,\\
\e^{-(p-1)} & \mbox{for}\ b>0
\end{array}
\right.
\quad
\mbox{if}\ \int_\R g(x)dx\neq0
\end{equation}
and
\begin{equation}
\label{lifespan_zeroKMT}
T(\e)\sim
\left\{
\begin{array}{ll}
\e^{-p(p-1)/(1-pb)} & \mbox{for}\ b<0,\\
\psi^{-1}(\e^{-p(p-1)}) & \mbox{for}\ b=0,\\
\e^{-p(p-1)} & \mbox{for}\ b>0
\end{array}
\right.
\quad
\mbox{if}\ \int_\R g(x)dx=0,
\end{equation}
where $\phi^{-1}$ and $\psi^{-1}$ are inverse functions defined by $\phi(s)=s\log(2+s)$ and $\psi(s)=s\log^p(2+s)$, respectively. 
The differences between $F=(1+t)^{-(p-1)}$ and $F=\LR{x}^{-(1+b)}$ in the lifespan estimates come from the decay property of the weights near the origin.
Indeed, for the time weights, there is a possibility that we get the global existence for $p>3$
due to the decay of the nonlinear term even if $|x|$ is small.
But there is no such a situation for the spatial weights.

Finally, we remark the setting of the weight function $F$ in our problem in (\ref{IVP}).
If $F=\LR{t+x}^{-(1+a)}$, we cannot expect to obtain the global existence result, because we have no decay property along $t+x=0$ when $x<0$. 
%For this reason, the result seems to be the same as that of the spatial weights. 
On the other hand, when $x>0$, we have the decay property 
for $\LR{t+x}^{-(1+a)}$ near the origin. For this reason, we get the lifespan estimates similar to the time weights.
Therefore, it is necessary to consider the weighted functions $\LR{t+|x|}^{-(1+a)}$ and $\LR{t-|x|}^{-(1+b)}$ separately. 
In addition, if $F \in C^1(\R \times (0,\infty))$ does not hold, we have no chance to obtain the classical solution of (\ref{IVP}) even locally in time.

Before stating our main results, we assume some conditions for the initial data.
Let $(f,g)\in C_{0}^2(\R)\times C_{0}^1(\R)$ satisfy 
\begin{equation}
\label{supp(f,g)}
{\rm supp}\ (f, g) \subset \{x\in \R : \ |x|\le R\} \ (R>1).
\end{equation}

Then, we obtain our main results. The global existence result (\ref{main_global}) is stated in Theorem \ref{thm_global1}. The lifespan estimates 
(\ref{main_lifespan1}) and (\ref{main_lifespan2}) are split into the following four theorems, Theorems \ref{thm_lower1}, \ref{thm_lower2}, \ref{thm_upper1} and \ref{thm_upper2}.
\begin{thm}\label{thm_global1}
Let $a+b>0$ and $a>0$. Assume the support condition (\ref{supp(f,g)}).
Then, there exists a positive constant $\e_1=\e_1(f,g,p,a,b,R)$ such that the problem (\ref{IVP}) 
admits a unique global-in-time solution $u\in C^2(\R\times[0,\infty))$ for $0<\e \le \e_1$.
\end{thm}
\begin{thm}\label{thm_lower1}
Assume the support condition (\ref{supp(f,g)}) and 
\begin{equation}
\label{initial_non-zero}
\int_{\R}g(x)dx\neq0.
\end{equation}
Then, there exists a positive constant $\e_2=\e_2(f,g,p,a,b,R)$ such that the problem (\ref{IVP}) admits a unique solution $u\in C^2(\R\times[0,T))$ as far as $T$ satisfies 
\[
T\le
\left\{
\begin{array}{lllll}
\exp(c\e^{-(p-1)}) &\mbox{if}&  a+b=0\ \mbox{and}\ a>0, \ \mbox{or}\ a=0\ \mbox{and} \ b>0, \\ 
\exp(c\e^{-(p-1)/2}) &\mbox{if}& a=b=0, \\
c\e^{-(p-1)/(-a)} &\mbox{if}& a<0 \ \mbox{and}\ b>0,\\
\phi_1^{-1}(c\e^{-(p-1)})&\mbox{if}& a<0\ \mbox{and}\ b=0,\\
c\e^{-(p-1)/(-a-b)}&\mbox{if}& a+b<0\ \mbox{and}\ b<0,\\
\end{array}
\right.
\]
where $0<\e \le \e_2$, $c$ is a positive constant independent of $\e$ and $\phi_1$ is the function in (\ref{phi1}).
\end{thm}

%\begin{thm}\label{thm_global2}
%Assume that (\ref{initial_zero}) holds. Let $a+b>0$ and $a>0$. 
%Then, there exists a positive constant $\e_3=\e_3(f,g,p,a,b,R)>0$ such that the problem (\ref{IVP}) 
%admits a unique solution $u\in C^2(\R\times[0,\infty))$ for $0<\e \.
%\end{thm}

\begin{thm}\label{thm_lower2}
Assume the support condition (\ref{supp(f,g)}) and 
\begin{equation}
\label{initial_zero}
\int_{\R}g(x)dx=0.
\end{equation}
Then, there exists a positive constant $\e_3=\e_3(f,g,p,a,b,R)$ such that the problem (\ref{IVP}) admits a unique solution 
$u\in C^2(\R\times[0,T))$ as far as $T$ satisfies 
\[
T\le
\left\{
\begin{array}{lllll}
\exp(c\e^{-(p-1)}) &\mbox{if}&  a=0\ \mbox{and}\ b>0,\\
\exp(c\e^{-p(p-1)}) &\mbox{if}& a+b=0\ \mbox{and} \ a>0, \\  
\exp(c\e^{-p(p-1)/(p+1)}) &\mbox{if}& a=b=0, \\
c\e^{-(p-1)/(-a)} &\mbox{if}& a<0 \ \mbox{and}\ b>0,\\
\psi_1^{-1}(c\e^{-p(p-1)})&\mbox{if}& a<0\ \mbox{and}\ b=0,\\
c\e^{-p(p-1)/(-pa-b)} &\mbox{if}& a<0\ \mbox{and}\ b<0,\\
\psi_2^{-1}(c\e^{-p(p-1)})&\mbox{if}& a=0\ \mbox{and}\ b<0,\\
c\e^{-p(p-1)/(-a-b)} &\mbox{if}& a+b<0 \ \mbox{and}\ a>0,
\end{array}
\right.
\]
where $0<\e \le \e_3$, $c$ is a positive constant independent of $\e$, $\psi_1$ and $\psi_2$ are defined in (\ref{psi1_psi2}).
\end{thm}

%kitamura_start%%%%%%%%%%%%%%%%%%%%%%%%%%%%%%%%%%%%%%%%%%%%%%%%%%%%%%%%%

\begin{thm}
\label{thm_upper1}
Assume the support condition (\ref{supp(f,g)})
and 
\begin{equation}
\label{positive_non-zero}
\int_{\R} g(x)dx > 0.
\end{equation}
Then, there exists a positive constant $\e_4=\e_4(f,g,p,a,b,R)$ such that
a classical solution $u\in C^2(\R\times[0,T))$ of (\ref{IVP}) cannot exist whenever $T$ satisfies
\[
T\ge
\left\{
\begin{array}{lllll}
\exp(C\e^{-(p-1)}) &\mbox{if}&  a+b=0\ \mbox{and}\ a>0 \ \mbox{or}\ a=0\ \mbox{and} \ b>0, \\ 
\exp(C\e^{-(p-1)/2}) &\mbox{if}& a=b=0, \\
C\e^{-(p-1)/(-a)} &\mbox{if}& a<0 \ \mbox{and}\ b>0,\\
\phi_1^{-1}(C\e^{-(p-1)})&\mbox{if}& a<0\ \mbox{and}\ b=0,\\
C\e^{-(p-1)/(-a-b)}&\mbox{if}& a+b<0\ \mbox{and}\ b<0,\\
\end{array}
\right.
\]
where $0<\e\le\e_4$, $C$ is a positive constant independent of $\e$
and $\phi_1$ is the function in (\ref{phi1}).
\end{thm}

\begin{thm}
\label{thm_upper2}
Assume the support condition (\ref{supp(f,g)})
and
\begin{equation}
\label{positive_zero}
f(x)\ge0(\not\equiv0),\quad g(x)\equiv0.
\end{equation}
Then, there exists a positive constant $\e_5=\e_5(f,p,a,b,R)$ such that
a classical solution $u\in C^2(\R\times[0,T))$ of (\ref{IVP}) cannot exist whenever $T$ satisfies
\[
T\ge
\left\{
\begin{array}{lllll}
\exp(C\e^{-(p-1)}) &\mbox{if}&  a=0\ \mbox{and}\ b>0,\\
\exp(C\e^{-p(p-1)}) &\mbox{if}& a+b=0\ \mbox{and} \ a>0, \\  
\exp(C\e^{-p(p-1)/(p+1)}) &\mbox{if}& a=b=0, \\
C\e^{-(p-1)/(-a)} &\mbox{if}& a<0 \ \mbox{and}\ b>0,\\
\psi_1^{-1}(C\e^{-p(p-1)})&\mbox{if}& a<0\ \mbox{and}\ b=0,\\
C\e^{-p(p-1)/(-pa-b)} &\mbox{if}& a<0\ \mbox{and}\ b<0,\\
\psi_2^{-1}(C\e^{-p(p-1)})&\mbox{if}& a=0\ \mbox{and}\ b<0,\\
C\e^{-p(p-1)/(-a-b)} &\mbox{if}& a+b<0 \ \mbox{and}\ a>0,
\end{array}
\right.
\]
where $0<\e\le\e_5$, $C$ is a positive constant independent of $\e$,
$\psi_1$ and $\psi_2$ are defined in (\ref{psi1_psi2}).
\end{thm}

Here, we address our lifespan estimates in the $(a,b)$-plane for the convenience of explanation. 
Figure 1 is related to Theorems \ref{thm_global1}, \ref{thm_lower1} and \ref{thm_upper1}.
Figure 2 is related to Theorems \ref{thm_global1}, \ref{thm_lower2} and \ref{thm_upper2}.

%\begin{tikzpicture} 
% \fill[white!40!lightgray](0.1,0)--(3,-2.9)--(3,3)--(0.1,3)--cycle;
% \fill[lightgray](-0.1,0.1)--(-3,0.1)--(-3,3)--(-0.1,3)--cycle;
% \fill[lightgray](0,-0.1)--(-3,-0.1)--(-3,-2.8)--(2.7,-2.8)--cycle;
% \draw[->,>=stealth,semithick] (-3.1,0)--(3.1,0)node[above]{$a$}; %x軸
% \draw[->,>=stealth,semithick] (0,-3.1)--(0,3.1)node[right]{$b$}; %y軸
% \draw (-3,3.22)node[right]{Figure 1.};
% \draw (3,3)node[below left]{$\infty$}; 
%  \draw[->] (0.2,0.2)--(0.03,0.03);
% \draw (0.1,0.3)node[right]{\rm exp$(C \e^{-\frac{p-1}{2}})$};
% %\draw (0,0)node[below left]{O}; %原点
% \draw[very thick,domain=0:3] plot(\x,-\x)node[left]{$a+b=0$};
%\end{tikzpicture}
%\begin{tikzpicture} 
% \fill[white!40!lightgray](0.1,0)--(3,-2.9)--(3,3)--(0.1,3)--cycle;
% \fill[lightgray](-0.1,0.1)--(-3,0.1)--(-3,3)--(-0.1,3)--cycle;
% \fill[lightgray](-0.1,-0.1)--(-3,-0.1)--(-3,-2.8)--(-0.1,-2.8)--cycle;
% \fill[lightgray](0.1,-0.2)--(0.1,-2.8)--(2.7,-2.8)--cycle;
% \draw[->,>=stealth,semithick] (-3.1,0)--(3.1,0)node[above]{$a$}; %x軸
%\draw[->,>=stealth,semithick] (0,-3.1)--(0,3.1)node[right]{$b$}; %y軸
% \draw (-3,3.22)node[right]{Figure 2.};
% \draw (3,3)node[below left]{$\infty$}; 
% \draw[->] (0.2,0.2)--(0.03,0.03);
% \draw (0.1,0.3)node[right]{\rm exp$(C\e^{-\frac{p(p-1)}{p+1}})$};
% %\draw (0,0)node[below left]{O}; %原点
% \draw[very thick,domain=0:3] plot(\x,-\x)node[left]{$a+b=0$};
%\end{tikzpicture}
\noindent
\begin{tikzpicture} 
 \fill[white!40!lightgray](0.1,0)--(2.9,-2.8)--(3.3,-2.8)--(3.3,3)--(0.1,3)--cycle;
 \fill[lightgray](-0.1,0.1)--(-3,0.1)--(-3,3)--(-0.1,3)--cycle;
 \fill[lightgray](0,-0.1)--(-3,-0.1)--(-3,-2.8)--(2.7,-2.8)--cycle;
 \draw[->,>=stealth,semithick] (-3.1,0)--(3.4,0)node[above]{$a$}; %x軸
 \draw[->,>=stealth,semithick] (0,-3.1)--(0,3.15)node[right]{$b$}; %y軸
 %\draw (-3,3.22)node[right]{Figure 1.};
 \fill [white] (1.9,0.98) rectangle (2.4,1.33);
 \draw (1.8,0.9)node[above right]{$\infty$}; 
 \draw[->] (0.67,-0.33)--(0.5,-0.5);
 \draw (0.6,-0.3)node[right]{\rm exp$(C\e^{-(p-1)})$};
 \draw[->] (0.2,0.2)--(0.03,0.03);
 \draw (0.1,0.3)node[right]{\rm exp$(C\e^{-\frac{p-1}{2}})$};
 \draw[->] (0.22,2)--(0,2);
 \draw (0.13,2)node[right]{\rm exp$(C\e^{-(p-1)})$};
 \fill [white] (-0.8,1.3) rectangle (-2.22,1.9);
 \draw (-2.3,1.2)node[above right]{$C\e^{-\frac{p-1}{-a}}$}; 
 \draw[->] (-2.5,0.26)--(-2.5,0);
 \draw (-2.77,0.5)node[right]{$\phi_{1}^{-1}(C\e^{-(p-1)})$};
 \fill [white] (-1.2,-1.7) rectangle (0.38,-1.1);
 \draw (-1.3,-1.8)node[above right]{$C\e^{-\frac{p-1}{-a-b}}$};
 \draw[very thick,domain=0:3] plot(\x,-\x)node[left]{$a+b=0$};
 \draw (0,-3.75)node[align=center]{\footnotesize Fig. 1. Lifespan estimates of \\ \footnotesize Theorems \ref{thm_global1}, \ref{thm_lower1} and \ref{thm_upper1}.};
\end{tikzpicture}
\begin{tikzpicture} 
 \fill[white!40!lightgray](0.1,0)--(2.9,-2.8)--(3.3,-2.8)--(3.3,3)--(0.1,3)--cycle;
 \fill[lightgray](-0.1,0.1)--(-3,0.1)--(-3,3)--(-0.1,3)--cycle;
 \fill[lightgray](-0.1,-0.1)--(-3,-0.1)--(-3,-2.8)--(-0.1,-2.8)--cycle;
 \fill[lightgray](0.1,-0.2)--(0.1,-2.8)--(2.7,-2.8)--cycle;
 \draw[->,>=stealth,semithick] (-3.1,0)--(3.4,0)node[above]{$a$}; %x軸
 \draw[->,>=stealth,semithick] (0,-3.1)--(0,3.15)node[right]{$b$}; %y軸
% \draw (-3,3.22)node[right]{Figure 2.};
 \fill [white] (1.9,0.98) rectangle (2.4,1.33);
 \draw (1.8,0.9)node[above right]{$\infty$}; 
 \draw[->] (0.67,-0.33)--(0.5,-0.5);
 \draw (0.55,-0.3)node[right]{\rm exp$(C\e^{-p(p-1)})$};
 \draw[->] (0.2,0.2)--(0.03,0.03);
 \draw (0.1,0.3)node[right]{\rm exp$(C\e^{-\frac{p(p-1)}{p+1}})$};
 \draw[->] (0.22,2)--(0,2);
 \draw (0.13,2)node[right]{\rm exp$(C\e^{-(p-1)})$};
 \fill [white] (-0.8,1.3) rectangle (-2.22,1.9);
 \draw (-2.3,1.2)node[above right]{$C\e^{-\frac{p-1}{-a}}$}; 
 \draw[->] (-2.67,0.26)--(-2.67,0);
 \draw (-2.94,0.5)node[right]{$\psi_1^{-1}(C\e^{-p(p-1)})$};
 \fill [white] (-0.5,-0.5) rectangle (-2.22,-1.1);
 \draw (-2.3,-1.2)node[above right]{$C\e^{-\frac{p(p-1)}{-pa-b}}$}; 
 \draw[->] (-0.22,-2)--(0,-2);
 \draw (-0.13,-2)node[left]{$\psi_2^{-1}(C\e^{-p(p-1)})$};
 \fill [white] (0.2,-2.5) rectangle (1.85,-1.92);
 \draw (0.1,-2.6)node[above right]{$C\e^{-\frac{p(p-1)}{-a-b}}$}; 
 \draw[very thick,domain=0:3] plot(\x,-\x)node[left]{$a+b=0$};
 \draw (0,-3.75)node[align=center]{\footnotesize Fig. 2. Lifespan estimates of \\ \footnotesize Theorems \ref{thm_global1}, \ref{thm_lower2} and \ref{thm_upper2}.};
\end{tikzpicture}

\begin{rem}
\label{rem2}
We  compare all the lifespan estimates when $a=-1$ in our results as for $F=\LR{t-\LR{x}}^{-(1+b)}$with those for $F=\LR{x}^{-(1+b)}$ by Kitamura, Morisawa and Takamura\cite{KMT}. If the total integral of $g$ does not vanish, 
they coincide with each other. However, if the total integral of $g$ vanishes, it does not hold. 
The reason for this situation can be described as follows.
The triangle domain which has a vertex $(x,t)$ in the following figure is the domain of the integral of the 
Duhamel term, (\ref{nonlinear}) below.  

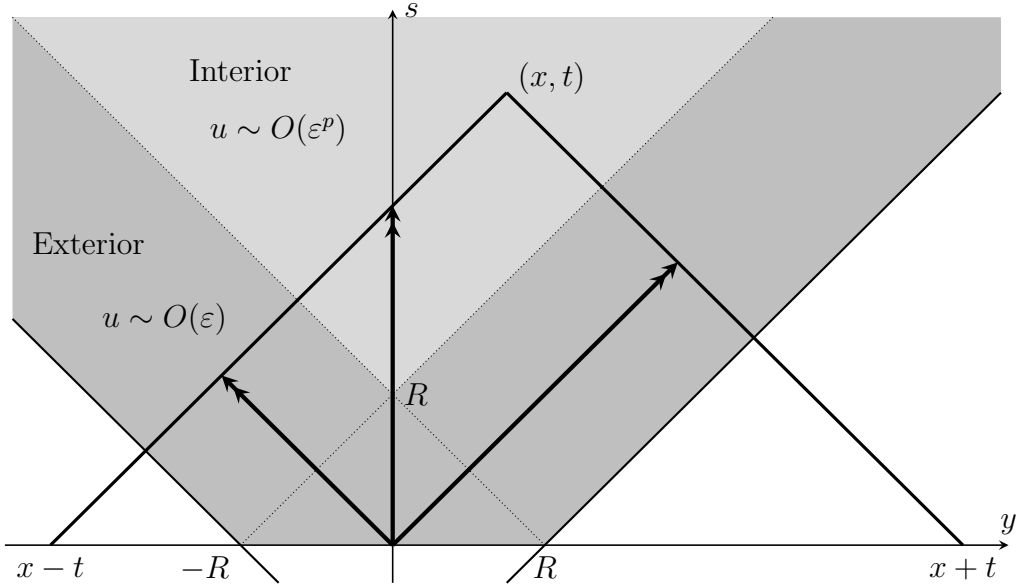
\noindent
\begin{tikzpicture} 
 \fill[lightgray](-2,0)--(2,0)--(8,6)--(8,7)--(-5,7)--(-5,3)--cycle;
 \fill[white!40!lightgray](0,2)--(-5,7)--(5,7)--cycle;
 %\fill[lightgray](0,-0.1)--(-3,-0.1)--(-3,-2.8)--(2.7,-2.8)--cycle;
 \draw[->,>=stealth,semithick] (-5.1,0)--(8.1,0)node[above]{$y$}; %x軸
 \draw[->,>=stealth,semithick] (0,-0.5)--(0,7.1)node[right]{$s$}; %y軸
 \draw (2,0)node[below]{$R$};
 \draw (-2,0)node[below left]{$-R$};
 \draw[thick,domain=-5:-1.5] plot(\x,-\x-2);
 \draw[thick,domain=1.5:8] plot(\x,\x-2);
 \draw[very thick,domain=1.5:7.5] plot(\x,-\x+7.5)node[below]{$x+t$};
 \draw[very thick,domain=1.5:-4.5] plot(\x,\x+4.5)node[below]{$x-t$};
 \draw[densely dotted,domain=5:-2] plot(\x,\x+2);
 \draw[densely dotted,domain=-5:2] plot(\x,-\x+2);
 \draw (0,2)node[right]{$R$};
 \draw (1.5,6.2)node[right]{$(x,t)$};
 \draw (-5,7)node[above right]{};%The case of $\d \int_\R g(y) dy =0$};
 \draw (-1.5,5.5)node{$u\sim O(\e^p)$};
 \draw (-4,4)node{{\rm Exterior}};
 \draw (-2,6.3)node{{\rm Interior}};
 \draw (-3,3)node{$u\sim O(\e)$};
 \draw [->>,>=stealth,ultra thick] (0,0)--(0,4.5);
\draw [->>,>=stealth,ultra thick] (0,0)--(3.75,3.75);
 \draw [->>,>=stealth,ultra thick] (0,0)--(-2.25,2.25);
\end{tikzpicture}
\par\noindent
When $b>-1$, we cannot derive any decay of the nonlinear term from the spatial weights $\LR{y}^{-(b+1)}$ along $y=0$. 
On the other hand, we cannot derive any decay of the nonlinear term from characteristic weights $\LR{s-|y|}^{-(b+1)}$ along $s-|y|=0$.
When $b<-1$, for the weights $\LR{y}^{-(b+1)}$, the growth of the solution appears on $s-|y|=0$. On the other hands, for the weights $\LR{s-|y|}^{-(b+1)}$, the growth of the solution appears on $y=0$. Keeping this situation in mind, we shall step into our purpose.

Set $a=-1$ as $F=\LR{t-|x|}^{-(1+b)}$. If we assume (\ref{initial_non-zero}), then the Huygens' principle (Lemma \ref{lm:huygens}) does not hold. 
Hence, we have that $u\sim O(\e)$ in its support. See the construction of the solution in the function space $Y_1$ below at the end of section 3. This fact helps us understand that the lifespan estimates in Theorem \ref{thm_lower1} and Theorem \ref{thm_upper1} 
are the same as that of (\ref{lifespan_non-zeroKMT}). 

In contrast, if we assume (\ref{initial_zero}) which implies that Huygens' principle holds, we have that $u \sim O(\e^p)$ in the interior domain $\{t+|x|\ge R$ and $t-|x|\ge R\}$ while we have that $u \sim O(\e)$ in the exterior domain $\{t+|x|\ge R$ and $|t-|x||\le R\}$. See the construction of the solution  in the function space $Y_2$ below at the end of section 4.
For $F=\LR{x}^{-(b+1)}$ $(b>-1)$, the solution $u$ does not decay along $y=0$ which is located in the interior domain. On the other hands, for $F=\LR{t-|x|}^{-(b+1)}$, the solution $u$ does not decay along $s-|y|=0$ which is located in the exterior domain. This fact helps us understand that the lifespan estimates in Theorems \ref{thm_lower2} and \ref{thm_upper2} are smaller than those of (\ref{lifespan_zeroKMT}) when $b>-1$. 
Conversely, in view of the considerations above, when $b<-1$, the lifespan estimates in Theorems \ref{thm_lower2} and \ref{thm_upper2} are larger than those of (\ref{lifespan_zeroKMT}).
%When $b=-1$, the lifespan estimates in our theorems are the same as that of (\ref{Zhou-est}).
\end{rem}
\begin{rem}
We remark that the weight function $F$ of $t+x$ and $t-x$ can cover a large class of power-type weights. 
For example, let us consider $F=\LR{t+x^2}^{-(a+1)}$. When $a>-1$, we may have the same results as for $F=\LR{t}^{-(a+1)}$ because $x^2$ in $F$ has a minor effect.
%Moreover, for the case $-1<a\le 1$, the lifespan estimates are the same as that of the time weights $\LR{t}^{-(a+1)}$. 
On the other hands, when $a<-1$, we may have the same results as for $F=\LR{x^2}^{-(1+a)}$ because $x^2$ in $F$ has a major effect.
It is sufficient to consider the linear combinations of $x$ and $t$ in the weighted functions.
%Thus, it is necessary that the variables in the weighted functions are the first order with respect to $x$ and $t$ because of the scaling for the wave equations.

Let us consider $F=\LR{t-C|x|}^{-(b+1)}$, where $C$ is a positive constant. 
If $0<C<1$ and $t-|x|>0$, then we have that 
\[
t\ge t-C|x|=(1-C)t+C(t-|x|)\ge (1-C)t,
\] so that the lifespan estimates may coincide with that for the time weights.  
On the other hands, if $C>1$, the line $s-C|y|=0$ in $(s,y)$-domain of the integral of the Duhamel term is located in the interior domain in Remark 1.1. Therefore, the lifespan estimates may coincide with that for the spatial weights.
%Therefore, it is enough to consider the weighted functions which have $x$, $t$ and $t-x$ variables to construct the general theory.
\end{rem}

This work was almost completed during the period when the first author was in the master course of Mathematical Institute of Tohoku University and the second author had the second affiliation of
Research Alliance Center of Mathematical Sciences, Tohoku University.
This paper is organized as follows. In the next section, we set the notation and prove some lemmas needed later. 
The proofs of Theorem \ref{thm_global1} with (\ref{initial_non-zero}) and Theorem \ref{thm_lower1} are established in Section 3. In Section 4, we prove Theorem \ref{thm_global1} 
with (\ref{initial_zero})  and Theorem \ref{thm_lower2}. Finaliy, we prove Theorem \ref{thm_upper1} and Theorem \ref{thm_upper2} in Section 5. 
All the proofs in this paper are based on the point-wise estimates of solutions which are originally introduced by John \cite{J79}.

%kitamura_end%%%%%%%%%%%%%%%%%%%%%%%%%%%%%%%%%%%%%%%%%%%%%%%%%%%%%%%%%

\section{Preliminaries}
\quad \ \ In this section, we set the notation and prove some useful lemmas.
Assume that $u\in C^2(\R\times[0,T])$ is a solution to the Cauchy problem (\ref{IVP}). Then the following finite propagation speeds 
holds:
\begin{equation}
\label{support_sol}
{\rm supp}\ u\subset \{(x,t)\in \R \times[0,T]: \ |x|\le t+R\}.
\end{equation}
For the proof, see Appendix of John \cite{J90}.

We define 
\begin{equation}
\label{linear}
u^0(x,t):=\frac{1}{2}\{f(x+t)+f(x-t)\}+\frac{1}{2}\int_{x-t}^{x+t}g(y)dy
\end{equation}
and
\begin{equation}
\label{nonlinear}
L(V)(x,t):=\frac{1}{2}\int_{0}^{t}\int_{x-t+s}^{x+t-s}\frac{V(y,s)}{\LR{s+\LR{y}}^{1+a}\LR{s-\LR{y}}^{1+b}}dyds
\end{equation}
for $V\in C(\R\times[0,T])$.

Let $(f,g)\in C^2(\R)\times C^1(\R)$. If $u\in C(\R\times[0,T])$ is a solution of 
\begin{equation}
\label{integral}
u(x,t)=\e u^0(x,t)+L(|u|^p)(x,t) \quad (x,t)\in \R\times[0,T],
\end{equation}
then $u\in C^2(\R\times[0,T])$ is the solution to the Cauchy problem (\ref{IVP}).

The following lemma is so-called Huygens' principle which is essential for the proof of main theorems.
\begin{lem}
\label{lm:huygens}
Assume that (\ref{supp(f,g)}) and (\ref{initial_zero}) hold. Then $u^0(x,t)$ satisfies
\begin{equation}
\label{huygens}
{\rm supp}\ u^0\subset \{(x,t)\in \R \times[0,\infty): \ (t-R)_{+}\le |x|\le t+R\},
\end{equation}
where $(a)_{+}=\max (0,a)$ for $a\in \R$. 
\end{lem}
For the proof, see Proposition 2.2 in \cite{KMT}. 

%kitamura_start%%%%%%%%%%%%%%%%%%%%%%%%%%%%%%%%%%%%%%%%%%%%%%%%%%%%%
\begin{lem}
\label{weight_equivalent}
For any $(x,t) \in \R \times [0,\infty)$, we have
%\begin{equation}
\[
%\label{equivalent1}
3^{-1}(1+|t-|x||) \leq \LR{t-\LR{x}} \leq \sqrt{2}  (1+|t-|x||) 
%\end{equation}
\]
and
%\begin{equation}
%\label{equivalent2}
\[
2^{-1}(1+t+|x|) \leq  \LR{t+\LR{x}} \leq \sqrt{2} (1+t+|x|).
\]
%\end{equation}
\end{lem}
\par\noindent
{\bf Proof.}
\noindent
We only prove $3^{-1}(1+|t-|x||) \leq \LR{t-\LR{x}}$ in the first inequalities. 
Others are trivial by taking the difference between the squares of right-hand side and
left-hand side. It is easy to see that
\[
\begin{aligned}
& 9\{1+(t-\sqrt{1+x^2})^2\} - (1+|t-|x||)^2 \\
& = 9+ 9t^2 -18t\sqrt{1+x^2} + 9(1+x^2) - 1 -2|t-|x||-t^2 + 2t|x| -x^2 \\
& = 9+ 7t^2 -18t\sqrt{1+x^2} +7(1+x^2) + (1-|t-|x||)^2 +4t|x| \\
& = h(x,t) + (1-|t-|x||)^2,
\end{aligned}
\]
where we set
\[
h(x,t):=7t^2 + 2(2|x|-9\sqrt{1+x^2})t + \{9+7(1+x^2)\}.
\]
Then, we see
\[
h(x,t)=7\left\{t+\frac{(2|x|-9\sqrt{1+x^2})}{7}\right\}^2+\frac{36|x|\sqrt{1+x^2}-36x^2+31}{7}>0.
\]
\if0
First, we define 
%\[
%h(x,t) := 9+7t^2+4t|x|+7(1+x^2)-18t\sqrt{1+x^2}.
%\]
%To summarize about $t$,
\[
h(x,t)= 7t^2 + 2(2|x|-9\sqrt{1+x^2})t + \{9+7(1+x^2)\}. 
\]
The discriminant for $h(t,x) = 0$ is as follows:
\[
\begin{aligned}
D/4 
& = (2|x|-9\sqrt{1+x^2})^2-7(7(1+x^2)+9) \\
& = 4x^2 -36|x|\sqrt{1+x^2} + 81(1+x^2)-49(1+x^2)-63 \\
& = <0 
\end{aligned}
\]
for $x \in \R$. This implies that $h(x,t) > 0$ for any $(x,t) \in \R \times [0,\infty)$, and
we finally obtain 
\[
\begin{aligned}
& 9(1+(t-\sqrt{1+x^2})^2) - (1+|t-|x||)^2 \\
& = 9+ 9t^2 -18t\sqrt{1+x^2} + 9(1+x^2) - 1 -2|t-|x||-t^2 + 2t|x| -x^2 \\
& = 9+ 7t^2 -18t\sqrt{1+x^2} +7(1+x^2) + (1-|t-|x||)^2 +4t|x| \\
& = h(x,t) + (1-|t-|x||)^2 > 0.
\end{aligned}
\]
\fi
\hfill
$\Box$
%kitamura_end%%%%%%%%%%%%%%%%%%%%%%%%%%%%%%%%%%%%%%%%%%%%%%%%%%%%%%%%%%%%%%%%%%%%%%%%%%%%%%%%%%%%%

Finally, we introduce the following domains to estimate the solutions:
\par\noindent
\begin{eqnarray*}
%\label{Int}
D_{\rm Int} :=& \{(x,t)\in \R \times[0,T]: \ t+|x|\ge R,\   t-|x|\ge R\},\\
%\label{ext}
D_{\rm Ext} :=& \{(x,t)\in \R \times[0,T]: \ t+|x|\ge R,\  |t-|x||\le R\},\\
%\label{ori}
D_{\rm Ori} :=& \{(x,t)\in \R \times[0,T]: \ t+|x|\le R,\   |t-|x||\le R\}.
\end{eqnarray*}
Note that the lifespan is determined by the point-wise estimates of the solution in $D_{\rm Int}$. 

\section{Proofs of Theorem \ref{thm_global1} and Theorem \ref{thm_lower1}}

We define $L^{\infty}$-norm of $u$ by
\begin{equation}
\label{norm1}
\|u\|_1:=\sup_{(x,t)\in \R\times[0,T]}|u(x,t)|.
\end{equation}

The following a priori estimate plays a key role in the proofs of Theorem \ref{thm_global1} and Theorem \ref{thm_lower1}.
\begin{lem}
\label{lem:apriori0}
Let $L$ be the linear integral operator defined by (\ref{nonlinear}).
Suppose that the assumptions of Theorem \ref{thm_global1} and Theorem \ref{thm_lower1} are fulfilled.
Assume that $u\in C(\R\times[0,T])$ and ${\rm supp}\ u\subset \{|x|\le t+R\}$ holds. 
Then, there exists a positive constant $C_1=C_1(p,a,b,R)$ such that
\begin{equation}
\label{apriori0}
\|L(|u|^p)\|_1\le C_1\|u\|_1^pE_1(T),
\end{equation}
where $E_1(T)$ is defined by
\begin{equation}
\label{E_1}
E_1(T):=
\left\{
\begin{array}{lll}
1 &\mbox{if}& a+b>0\ \mbox{and}\ a>0,\\
\log(T+3R) &\mbox{if}& a+b=0\ \mbox{and}\ a>0, \\
& &  \mbox{or}\ a=0\ \mbox{and} \ b>0, \\ 
\log^2(T+3R) &\mbox{if}& a=b=0, \\
(T+2R)^{-a} &\mbox{if}& a<0 \ \mbox{and}\ b>0,\\
(T+R)^{-a}\log(T+3R)&\mbox{if}& a<0\ \mbox{and}\ b=0,\\
(T+2R)^{-(a+b)}&\mbox{if}& a+b<0\ \mbox{and}\ b<0.\\
\end{array}
\right.
\end{equation}

\end{lem}
\par\noindent
{\bf Proof.}
\noindent
We denote by $C$ a positive constant which depends only on $p,a,b,R$ and $j$ may change from line to line. 
The definition of (\ref{nonlinear}) gives us 
\[
\begin{array}{llll}
L(|u|^p)(x,t)&\le\d C\|u\|_1^p\int_{0}^{t}\int_{x-t+s}^{x+t-s}\frac{\chi_0(y,s)}{\LR{s+\LR{y}}^{1+a}\LR{s-\LR{y}}^{1+b}}dyds\\
&\le\d C \|u\|_1^p\int_{0}^{t}\int_{x-t+s}^{x+t-s}\frac{\chi_0(y,s)}{(1+|s+|y||)^{1+a}(1+|s-|y||)^{1+b}}dyds,
\end{array}
\]
where 
\begin{equation}
\label{chi}
\chi_0(y,s):=
\left\{
\begin{array}{ll}
1 & \mbox{if}\ |y|\le s+R,\\
0 & \mbox{otherwise}.
\end{array}
\right.
\end{equation}
Set
\begin{equation}
\label{I2}
I(x,t):=\int_{0}^{t}\int_{x-t+s}^{x+t-s}\frac{\chi_0(y,s)}{(1+|s+|y||)^{1+a}(1+|s-|y||)^{1+b}}dyds.
\end{equation}
Then, the inequality (\ref{apriori0}) follows from 
\begin{equation}
\label{basic0}
I(x,t)\le C_1 E_1(T)\quad \mbox{in} \ D_{\rm Int} \cup D_{\rm Ext}\cup D_{\rm Ori}. 
\end{equation}
Thus, we show (\ref{basic0}) in the following.  

Without loss of generality, we may assume $x\ge0$ due to the symmetry of $I(x,t)$ as $I(x,t) = I(-x, t)$ holds by (\ref{I2}). 
Changing variables by 
\begin{equation}
\label{coordinate}
\alpha=s+y, \ \beta=s-y
\end{equation}
in the integral of (\ref{I2}), we get
\[
I(x,t)\le\left\{
\begin{array}{lll}
I_{11}(x,t)+I_{12}(x,t)+I_{13}(x,t)+I_{14}(x,t) & \mbox{in} & D_{\rm Int}, \\
I_{21}(x,t)+I_{22}(x,t) & \mbox{in} & D_{\rm Ext}\cup D_{\rm Ori}.\\
%I_{31}(x,t)+I_{32}(x,t) &\mbox{in}& D_{ori}.
\end{array}
\right.
\]
Here we set
\[
\begin{array}{llll}
I_{11}(x,t)&\d:=\int_{R}^{t-x}\frac{d\beta}{(1+\beta)^{1+b}}\int_{\beta}^{t+x}\frac{d\alpha}{(1+\alpha)^{1+a}},\\
I_{12}(x,t)&\d:=\int_{R}^{t-x}\frac{d\alpha}{(1+\alpha)^{1+b}}\int_{\alpha}^{t-x}\frac{d\beta}{(1+\beta)^{1+a}},\\
I_{13}(x,t)&\d:=\int_{-R}^{R}\frac{d\beta}{(1+|\beta|)^{1+b}}\int_{|\beta|}^{t+x}\frac{d\alpha}{(1+\alpha)^{1+a}},\\
I_{14}(x,t)&\d:=\int_{-R}^{R}\frac{d\alpha}{(1+\alpha)^{1+b}}\int_{|\alpha|}^{t-x}
\frac{d\beta}{(1+\beta)^{1+a}},\\
I_{21}(x,t)&\d:=\int_{-R}^{t-x}\frac{d\beta}{(1+|\beta|)^{1+b}}\int_{|\beta|}^{t+x}\frac{d\alpha}{(1+\alpha)^{1+a}},\\
I_{22}(x,t)&\d:=\int_{-(t-x)_+}^{(t-x)_{+}}\frac{d\alpha}{(1+|\alpha|)^{1+b}}\int_{|\alpha|}^{t-x}\frac{d\beta}{(1+\beta)^{1+a}}.\\
%I_{31}(x,t)&\d=\int_{-R}^{t-x}\frac{d\beta}{(1+|\beta|)^{1+b}}\int_{\beta}^{t+x}\frac{d\alpha}{(1+|\alpha|)^{1+a}},\\
%I_{32}(x,t)&\d=\int_{0}^{t-x}\frac{d\beta}{(1+\beta)^{1+a}}\int_{-\beta}^{\beta}\frac{d\alpha}{(1+|\alpha|)^{1+b}}.\\
\end{array}
\]

We first state the following lemma which is useful to estimate the above integrals.
\begin{lem}
Let $q\in\R$ and $\mu, \nu\ge0$. Then there exists a positive constant $C=C(q)$ such that 
\begin{equation}
\label{est:alpha-beta}
\int_{\mu}^{\nu}\frac{d\xi}{(1+\xi)^{1+q}}\le
C\left\{
\begin{array}{ll}
(1+\mu)^{-q} &\mbox{if} \ q>0,\\ 
\log(\nu+1) &\mbox{if} \ q=0,\\
(\nu+1)^{-q} & \mbox{if} \ q<0.\\
\end{array}
\right.
\end{equation}
\end{lem}
The proof is easy so that we shall omit it.

\begin{lem}
Let $q<0$ and $R>1$. Assume that $(x,t)\in[0,\infty)\times[0,T]$ satisfies $x\le t+R$. 
Then there exists a positive constant $C=C(q,R)$ such that 
\begin{equation}
\label{est:T1}
\log(t+x+3R)\le C\log(T+3R)
\end{equation}
and
\begin{equation}
\label{est:T2}
(t+x+R)^{-q}\le C(T+R)^{-q}.
\end{equation}
\end{lem}
\noindent
{\bf Proof.}
We show the inequality (\ref{est:T1}) only as (\ref{est:T2}) is trivial. The assumptions give us 
\[
\log(t+x+3R)\le \log (2T+4R)\le 2\log(T+3R).
\]
\hfill
$\Box$

First,  we shall estimate $I_{11}$ in $D_{\rm Int}$.
\vskip10pt
\par\noindent
{\bf Case 1. $\v{a>0}$ and $\v{a+b>0}$, or $\v{a>0}$ and $\v{a+b=0}$, or $\v{a>0}$ and $\v{a+b<0}$.} 

Applying (\ref{est:alpha-beta}) with $q=a(>0)$, $\mu=\beta$, $\nu=t+x$ to the $\alpha$-integral, we have
\[
I_{11}(x,t)\le C\int_{R}^{t-x}\frac{d\beta}{(1+\beta)^{1+a+b}}.
\]
Since (\ref{est:alpha-beta}) with $q=a+b$, $\mu=R$ and $\nu=t-x$, we obtain
\[
I_{11}(x,t)\le 
C\left\{
\begin{array}{ll}
1 &\mbox{if} \ a+b>0,\\ 
\log(t-x+1) &\mbox{if} \ a+b=0,\\
(t-x+1)^{-(a+b)} & \mbox{if} \ a+b<0.\\
\end{array}
\right.
\] 
It follows from $t-x+1\le t+x+R$, (\ref{est:T1}) and (\ref{est:T2}) that 
\[
I_{11}(x,t)\le 
C\left\{
\begin{array}{ll}
1 &\mbox{if} \ a+b>0,\\ 
\log(T+3R) &\mbox{if} \ a+b=0,\\
(T+2R)^{-(a+b)} & \mbox{if} \ a+b<0.\\
\end{array}
\right.
\] 
\vskip10pt
\par\noindent
{\bf Case 2. $\v{a=0}$ and $\v{b>0}$, or $\v{a=0}$ and $\v{b=0}$, or $\v{a=0}$ and $\v{b<0}$.} 

First, we consider the case of $a=0$ and $b>0$, or $a=0$ and $b=0$.
Applying (\ref{est:alpha-beta}) with $q=a(=0)$, $\mu=\beta$, $\nu=t+x$ to the $\alpha$-integral, we get
\[
I_{11}(x,t)\le C\log(t+x+1)\int_{R}^{t-x}\frac{d\beta}{(1+\beta)^{1+b}}.
\]
Making use of (\ref{est:T1}) and (\ref{est:alpha-beta}) with $\mu=R$, $\nu=t-x$ and $q=b(\ge0)$, we obtain 
\[
I_{11}(x,t)\le
C\left\{
\begin{array}{lll}
\log(T+3R) &\mbox{if} \ a=0 \ \mbox{and} \ b>0,\\
\log^2(T+3R) &\mbox{if} \ a=0\ \mbox{and}\ b=0.
\end{array}
\right.
\]
Next, we consider the case of $a=0$ and $b<0$. Since $b<0$ and $\beta\le \alpha$, we have
\[
I_{11}(x,t)\le C \int_{R}^{t-x}\frac{d\beta}{(1+\beta)^{1+b/2}}\int_{\beta}^{t+x}\frac{d\alpha}{(1+\alpha)^{1+b/2}}.
\]
By virtue of (\ref{est:alpha-beta}) with $q=b/2 (<0)$ and (\ref{est:T2}), we get
\[
I_{11}(x,t)\le C(T+R)^{-b}.
\]
\vskip10pt
\par\noindent
{\bf Case 3. $\v{a<0}$ and $\v{b>0}$, or $\v{a<0}$ and $\v{b=0}$, or $\v{a<0}$ and $\v{b<0}$.} 

Making use of (\ref{est:alpha-beta}) with $q=a(<0)$ and (\ref{est:T2}) to the $\alpha$-integral, we obtain
\[
I_{11}(x,t)\le C (T+R)^{-a} \int_{R}^{t-x}\frac{d\beta}{(1+\beta)^{1+b}}.
\]
It follows from (\ref{est:alpha-beta}) with $q=b$ and (\ref{est:T2}) that
\[
I_{11}(x,t)\le
C\left\{
\begin{array}{lll}
(T+2R)^{-a} &\mbox{if} \ b>0,\\
\log(T+3R)(T+R)^{-a} &\mbox{if} \ b=0,\\ 
(T+2R)^{-(a+b)} & \mbox{if} \ b<0.
\end{array}
\right.
\]
Summing up all the estimates, we get
\[
I_{11}(x,t)\le CE_1(T)\quad \mbox{in}\ D_{\rm Int}.
\]

Next, we shall estimate $I_{12}$ in $D_{\rm Int}$. Noticing that $\beta\le t-x\le t+x$, the 
estimates of this integral are the same as the estimates for $I_{11}$, which is replaced by $\alpha$ and $\beta$. So, we have 
\[
I_{12}(x,t)\le CE_1(T)\quad \mbox{in}\ D_{\rm Int}.
\]

Next, we shall estimate $I_{13}$ and $I_{14}$ in $D_{\rm Int}$. 
We note that both the $\beta$-integral in $I_{13}$ and $\alpha$-integral in $I_{14}$ are bounded by $C$. Thus, we get
\[
I_{13}(x,t)\le C\int_{0}^{t+x}\frac{d\alpha}{(1+\alpha)^{1+a}} \ \ \mbox{and}\ \ I_{14}(x,t)\le C\int_{0}^{t+x}\frac{d\beta}{(1+\beta)^{1+a}}.
\]
Making use of (\ref{est:alpha-beta}) with $q=a$, (\ref{est:T1}) and (\ref{est:T2}) to the above integrals, we obtain 
\begin{equation}
\label{I34}
I_{13}(x,t), I_{14}(x,t)\le 
C\left\{
\begin{array}{ll}
1 &\mbox{if} \ a>0,\\ 
\log(T+3R) &\mbox{if} \ a=0,\\
(T+R)^{-a} & \mbox{if} \ a<0.\\
\end{array}
\right.
\end{equation}
It follows from
\begin{equation}
\label{TR-est}
T+R\ge1 \quad \mbox{and}\quad  \log(T+3R)\ge 1,
\end{equation}
(\ref{I34}) and the definition of $E_1(T)$ in (\ref{E_1}), that
\[
I_{13}(x,t),I_{14}(x,t)\le CE_1(T)\quad \mbox{in}\ D_{\rm Int}.
\]

Next, we shall estimate $I_{21}$ and $I_{22}$ in $D_{\rm Ext}\cup D_{\rm Ori}$. First, we investigate them in $D_{\rm Ext}$.
Since $|t-x|\le R$ holds for $(x,t)\in D_{\rm Ext}$, we get
\[
I_{21}(x,t)\le\int_{-R}^{R}\frac{d\beta}{(1+|\beta|)^{1+b}}\int_{\beta}^{t+x}\frac{d\alpha}{(1+|\alpha|)^{1+a}}.\\
\]
For $t-x>0$, we have
\[
I_{22}(x,t)\le\int_{-R}^{R}\frac{d\alpha}{(1+|\alpha|)^{1+b}}\int_{-R}^{R}\frac{d\beta}{(1+|\beta|)^{1+a}}\le  C.
\]
Making use of (\ref{TR-est}) for the estimates of $I_{22}$, we get
\[
I_{22}(x,t)\le CE_1(T)\quad \mbox{in}\ D_{\rm Ext}.
\]
For the $\alpha$-integral of $I_{21}$, we obtain 
\begin{equation}
\label{ext-est1}
\d \int_{-R}^{t+x}\frac{d\alpha}{(1+|\alpha|)^{1+a}}\le 
C\left\{
\begin{array}{lll}
1 &\mbox{if} \ a>0,\\ 
\log(t+x+3R) &\mbox{if} \ a=0,\\
(t+x+R)^{-a} & \mbox{if} \ a<0.
\end{array}
\right.
\end{equation}
where we have used $\beta\ge-R$.
It follows from (\ref{est:T1}), (\ref{est:T2}) and (\ref{TR-est}) that
\[
I_{21}(x,t)\le CE_1(T)\quad \mbox{in}\ D_{\rm Ext}.
\]
Finally, we shall estimate $I_{21}$ and $I_{22}$ in $D_{\rm Ori}$. 
Noticing that $|t-x|\le R$ and $t+x\le R$ hold for $(x,t)\in D_{\rm Ori}$, we obtain
\[
I_{21}(x,t)\le \int_{-R}^{R}\frac{d\beta}{(1+|\beta|)^{1+b}}\int_{-R}^{R}\frac{d\alpha}{(1+|\alpha|)^{1+a}}\le C
\]
and
\[
I_{22}(x,t)\le \int_{-R}^{R}\frac{d\alpha}{(1+|\alpha|)^{1+b}}\int_{-R}^{R}\frac{d\beta}{(1+|\beta|)^{1+a}}\le C
\]
for $t-x>0$.
Making use of (\ref{TR-est}), we get
\[
I_{21}(x,t), I_{22}(x,t)\le CE_1(T)\quad \mbox{in}\ D_{\rm Ori}.
\]
Summing up all the estimates, we obtain (\ref{basic0}). Therefore the proof of (\ref{apriori0}) is established.
\hfill
$\Box$
\vskip10pt
\par\noindent
{\bf Proofs of Theorem \ref{thm_global1} and Theorem \ref{thm_lower1}.}
\noindent
Let us define a Banach space
\[
X_1:=\{u\in C(\R\times[0,T]): \ {\rm supp}\ u\subset \{|x|\le t+R\} \}
\]
which is equipped with the norm (\ref{norm1}). 
Define a sequence of functions $\{u_n\}_{n\in \N}$ by 
\[
u_{n}=u_0+L(|u_{n-1}|^p),\quad u_0=\e u^0,
\]
where $L$ and $u^0$ are defined in (\ref{nonlinear}) and (\ref{linear}) respectively.

Set
\begin{equation}
\label{linear-C_0}
C_0:=\|f\|_{L^{\infty}(\R)}+\frac{1}{2}\|g\|_{L^1(\R)}>0. 
\end{equation}
We shall construct a solution of the integral equation (\ref{integral}) in a closed subspace 
\[
Y_1:=\{u\in X_1: \|u\|_1\le 2C_0\e\}
\] 
of $X_1$.
Note that $u_0\in Y_1$ holds. Because it follows from (\ref{linear}) and the assumption for $(f,g)$ that $|u_0(x,t)|\le \e C_0$. 
Then, analogously to the proof of Theorem 1.2 in \cite{W17} with $M:=C_0$, we can show that $u_n\in Y_1$ $(n\in \N)$ provided 
\[
2^pC_0^{p-1}C_1E_1(T)\e^{p-1}\le 1
\]
holds, where $C_1$ and $E_1(T)$ are defined in (\ref{apriori0}). 
Moreover, $\{u_n\}$ is a Cauchy sequence in $Y_1$ provided 
\[
2^{p+1}pC_0^{p-1}C_1E_1(T)\e^{p-1}\le 1
\]
holds. Since $Y_1$ is complete, 
there exists a function $u$ such that $u_n$ converges to $u$ in $Y_1$. Therefore $u$ satisfies (\ref{integral}).
\hfill
$\Box$

\section{Proofs of Theorem \ref{thm_global1} and Theorem \ref{thm_lower2}}

We define the following weight function:
\begin{equation}
\label{weight}
w(|x|,t):=\left\{
\begin{array}{llll}
1 &\mbox{if}& \ a>0,\\
\log^{-1}(t+|x|+3R) &\mbox{if}& \ a=0,\\
(t+|x|+3R)^{a} &\mbox{if}&\ a<0.\\
\end{array}
\right.
\end{equation}
Denote a weighted $L^{\infty}$-norm of $U$ by
\begin{equation}
\label{norm2}
\|U\|_2:=\sup_{(x,t)\in \R\times[0,T]}w(|x|,t)|U(x,t)|.
\end{equation}

Then we have the following a priori estimates.
\begin{lem}
\label{lem:apriori1}
Suppose that the assumptions of Theorem \ref{thm_global1} and Theorem \ref{thm_lower2} are fulfilled.
Assume that $U, U^0\in C(\R\times[0,T])$ with ${\rm supp}\ U\subset \{(x,t)\in \R \times[0,T]: \ |x|\le t+R\}$ and
 ${\rm supp}\ U^0\subset \{(x,t)\in \R \times[0,T]: \ (t-R)_{+}\le |x|\le t+R\}$ hold. 
Then, there exists a positive constant $C_2=C_2(p,a,b,R,j)$ such that
\begin{equation}
\label{apriori1}
\|L(|U^0|^{p-j}|U|^j)\|_2\le C_2(\|U\|_2D(T))^j \quad \mbox{for}\ j=0,1,
\end{equation}
where $D(T)$ is defined by
\[
D(T):=
\left\{
\begin{array}{lll}
1 & \mbox{if}\ a>0,\\
\log(T+3R) & \mbox{if}\ a=0,\\
(T+2R)^{-a} & \mbox{if}\ a<0.
\end{array}
\right.
\]
\end{lem}
\begin{lem}
\label{lem:apriori2}
Suppose that the assumptions of Theorem \ref{thm_global1} and Theorem \ref{thm_lower2} are fulfilled.
Assume that $U\in C(\R\times[0,T])$ with ${\rm supp}\ U\subset \{(x,t)\in \R \times[0,T]: |x|\le t+R\}$ holds. 
Then, there exists a positive constant $C_3=C_3(p,a,b,R)$ such that
\begin{equation}
\label{apriori2}
\|L(|U|^p)\|_2\le C_3\|U\|_2^pE_2(T),
\end{equation}
where $E_2(T)$ is defined by
\[
E_2(T):=
\left\{
\begin{array}{lllll}
1 &\mbox{if}&  a+b>0\ \mbox{and}\ a>0,\\
\log^p(T+3R) &\mbox{if}&  a=0\ \mbox{and}\ b>0,\\
\log(T+3R) &\mbox{if}& a+b=0\ \mbox{and} \ a>0, \\  
\log^{p+1}(T+3R) &\mbox{if}& a=b=0, \\
(T+2R)^{-pa} &\mbox{if}& a<0 \ \mbox{and}\ b>0,\\
(T+2R)^{-pa}\log(T+3R)&\mbox{if}& a<0\ \mbox{and}\ b=0,\\
(T+2R)^{-(pa+b)} &\mbox{if}& a<0\ \mbox{and}\ b<0,\\
\log^{p-1}(T+3R)(T+R)^{-b}&\mbox{if}& a=0\ \mbox{and}\ b<0,\\
(T+2R)^{-(a+b)} &\mbox{if}& a+b<0 \ \mbox{and}\ a>0.
\end{array}
\right.
\]
\end{lem}
\par\noindent
{\bf Proof of Lemma \ref{lem:apriori1}.}
\noindent
We denote by $C$ a positive constant which depends only on $p,a,b,R$ and $j$ may change from line to line. 
Making use of (\ref{huygens}), we have
\[
\begin{array}{llll}
&&\d L(|U^0|^{p-j}|U|^j)(x,t)\\
&&\le\d CC_0^{p-j}\|U\|_2^{j}\int_{0}^{t}\int_{x-t+s}^{x+t-s}\frac{w^{-j}(|y|,s)\chi_1(y,s)}{\LR{s+\LR{y}}^{1+a}
\LR{s-\LR{y}}^{1+b}}dyds\\
&&\le\d CC_0^{p-j}\|U\|_2^j\int_{0}^{t}\int_{x-t+s}^{x+t-s}\frac{w^{-j}(|y|,s)\chi_1(y,s)}{(1+|s+|y||)^{1+a}
(1+|s-|y||)^{1+b}}dyds,
\end{array}
\]
where 
\[
\chi_1(y,s):=
\left\{
\begin{array}{ll}
1 & \mbox{if}\ (s-R)_+\le |y|\le s+R,\\
0 & \mbox{otherwise},
\end{array}
\right.
\]
and $C_0$ is the constant in (\ref{linear-C_0}).

Set 
\[
J(x,t):=\int_{0}^{t}\int_{x-t+s}^{x+t-s}\frac{w^{-j}(|y|,s)\chi_1(y,s)}{(1+|s+|y||)^{1+a}(1+|s-|y||)^{1+b}}dyds.
\]
Then, the inequality (\ref{apriori1}) follows from 
%\begin{equation}
%\label{basic1}
\[
J(x,t)\le C_2w^{-1}(|x|,t) D^j(T)\quad \mbox{in} \ D_{\rm Int}\cup D_{\rm Ext}\cup D_{\rm Ori}. 
\]
%\end{equation}
Thus, we show the above inequality in the following.  

Without loss of generality, we may assume $x\ge0$ due to the symmetry of $J(x,t)$ in $x$. 
Changing variables by (\ref{coordinate}), we obtain
\[
J(x,t)\le
C\left\{
\begin{array}{llll}
J_{11}(x,t)+J_{12}(x,t)\  \mbox{in} \ D_{\rm Int}\\
J_{21}(x,t)+J_{22}(x,t)\  \mbox{in} \  D_{\rm Ext}\cup D_{\rm Ori}.
\end{array}
\right.
\]
Here we set
\[
\begin{array}{llll}
J_{11}(x,t)&\d:=\int_{-R}^{R}\frac{d\beta}{(1+|\beta|)^{1+b}}\int_{\beta}^{t+x}\frac{w^{-j}(|y|,s)}{(1+|\alpha|)^{1+a}}d\alpha,\\
J_{12}(x,t)&\d:=\int_{-R}^{R}\frac{d\alpha}{(1+|\alpha|)^{1+b}}\int_{|\alpha|}^{t-x}\frac{w^{-j}(-y,s)}{(1+\beta)^{1+a}}d\beta,\\
J_{21}(x,t)&\d:=\int_{-R}^{t-x}\frac{d\beta}{(1+|\beta|)^{1+b}}\int_{\beta}^{t+x}\frac{w^{-j}(|y|,s)}{(1+|\alpha|)^{1+a}}d\alpha,\\
J_{22}(x,t)&\d:=\int_{-(t-x)}^{(t-x)_{+}}\frac{d\alpha}{(1+|\alpha|)^{1+b}}\int_{|\alpha|}^{t-x}\frac{w^{-j}(-y,s)}{(1+\beta)^{1+a}}d\beta.\\
\end{array}
\]
First, we shall estimate $J_{11}$ in $D_{\rm Int}$. 
Noticing that
\[
w^{-1}(y,s)\le C\left\{
\begin{array}{llll}
1 &\mbox{if}& \ a>0,\\
\log(T+3R) &\mbox{if}& \ a=0,\\
(T+R)^{-a} &\mbox{if}&\ a<0\\
\end{array}
\right.
\]
by $\alpha\le t+x$, (\ref{est:T1}) and (\ref{est:T2}), we obtain
\[
J_{11}(x,t)\le CD^{j}(T) \int_{-R}^{R}\frac{d\beta}{(1+|\beta|)^{1+b}}\int_{-R}^{t+x}\frac{d\alpha}{(1+|\alpha|)^{1+a}}.
\]
Due to (\ref{ext-est1}), we get
\[
\int_{-R}^{t+x}\frac{d\alpha}{(1+|\alpha|)^{1+a}}\le w^{-1}(x,t)
\]
which implies that
\[
J_{11}(x,t)\le CD^{j}(T)w^{-1}(x,t).
\]
Next, we shall estimates $J_{12}$ in $D_{\rm Int}$. 
Because of $\beta\le t-x $, we get 
\[
w^{-1}(-y,s)\le C\left\{
\begin{array}{llll}
1 &\mbox{if}& \ a>0,\\
\log(T+3R) &\mbox{if}& \ a=0,\\
(T+R)^{-a} &\mbox{if}&\ a<0,\\
\end{array}
\right.
\]
where we have used  (\ref{est:T1}) and (\ref{est:T2}). Hence we obtain 
\[
\begin{array}{lll}
J_{12}(x,t)&\le\d CD^j(T)\int_{-R}^{R}\frac{d\alpha}{(1+|\alpha|)^{1+b}}\int_{|\alpha|}^{t-x}\frac{d\beta}{(1+\beta)^{1+a}}\\
&\le\d CD^j(T)\int_{-R}^{R}\frac{d\alpha}{(1+|\alpha|)^{1+b}}\int_{-R}^{t+x}\frac{d\beta}{(1+|\beta|)^{1+a}}.
\end{array}
\]
Making use of (\ref{ext-est1}), we obtain 
\[
J_{12}(x,t)\le CD^j(T)w^{-1}(x,t).
\]
We next estimate $J_{21}$ and $J_{22}$ in $D_{\rm Ext}$. 
Because of $|t-x|\le R$, we get 
\[
J_{21}(x,t)\le J_{11}(x,t),\quad J_{22}(x,t)\le J_{12}(x,t)\ \mbox{in} \ D_{\rm Ext}.
\]
Thus, the desired estimates for $J_{21}$ and $J_{22}$ are established. 
Finally, we shall estimate $J_{21}$ and $J_{22}$ in $D_{\rm Ori}$. 
Noticing that $|t-x|\le R$ and $t+x\le R$ hold for $(x,t)\in D_{\rm Ori}$, we obtain
\[
J_{1}(x,t)\le \int_{-R}^{R}\frac{d\beta}{(1+|\beta|)^{1+b}}\int_{-R}^{R}\frac{w^{-j}(y,s)}{(1+|\alpha|)^{1+a}}d\alpha\le C
\]
and
\[
J_{2}(x,t)\le \int_{-R}^{R}\frac{d\alpha}{(1+|\alpha|)^{1+b}}\int_{-R}^{R}\frac{w^{-j}(-y,s)}{(1+|\beta|)^{1+a}}d\beta\le C
\]
for $t-x>0$. It follows from $w(|x|,t)\ge1$ and (\ref{TR-est}) that 
\[
J_{21}(x,t), J_{22}(x,t)\le CD^j(T)w^{-1}(x,t) \quad \mbox{in}\ D_{\rm Ori}.
\]
Therefore, we get (\ref{apriori1}).
\hfill
$\Box$
\vskip10pt
\par\noindent
{\bf Proof of Lemma \ref{lem:apriori2}.}
\noindent
We denote by $C$ a positive constant which depends only on $p,a,b,R$ and $j$ may change from line to line. 
The definition of (\ref{nonlinear}) gives us 
\[
\begin{array}{llll}
L(|U|^p)(x,t)
&\le\d C\|U\|_2^p\int_{0}^{t}\int_{x-t+s}^{x+t-s}\frac{w^{-p}(|y|,s)\chi_0(y,s)}{\LR{s+\LR{y}}^{1+a}\LR{s-\LR{y}}^{1+b}}dyds\\
&\le\d C \|U\|_2^p\int_{0}^{t}\int_{x-t+s}^{x+t-s}\frac{w^{-p}(|y|,s)\chi_0(y,s)}{(1+|s+|y||)^{1+a}
(1+|s-|y||)^{1+b}}dyds,
\end{array}
\]
where $\chi_0$ is the function in (\ref{chi}).
Set
\[
\wt{J}(x,t):=\int_{0}^{t}\int_{x-t+s}^{x+t-s}\frac{w^{-p}(|y|,s)\chi_0(y,s)}{(1+|s+|y||)^{1+a}
(1+|s-|y||)^{1+b}}dyds.
\]
Then, the inequaltiy (\ref{apriori2}) follows from 
\begin{equation}
\label{basic2}
\wt{J}(x,t)\le C_3 w^{-1}(|x|,t)E_2(T)\quad \mbox{in} \ D_{\rm Int} \cup D_{\rm Ext}\cup D_{\rm Ori}. 
\end{equation}
Thus, we show (\ref{basic2}) in the following.  

Without loss of generality, we may assume $x\ge0$ due to the symmetry of $\wt{J}(x,t)$ in $x$. 
Changing variables by (\ref{coordinate}), we get
\[
\wt{J}(x,t)\le\left\{
\begin{array}{lll}
\wt{J_{11}}(x,t)+\wt{J_{12}}(x,t)+\wt{J_{13}}(x,t)+\wt{J_{14}}(x,t) & \mbox{in} & D_{\rm Int}, \\
\wt{J_{21}}(x,t)+\wt{J_{22}}(x,t) & \mbox{in} & D_{\rm Ext}\cup D_{\rm Ori}.\\
\end{array}
\right.
\]
Here we set
\[
\begin{array}{llll}
\wt{J_{11}}(x,t)&\d:=\int_{R}^{t-x}\frac{d\beta}{(1+\beta)^{1+b}}\int_{\beta}^{t+x}\frac{w^{-p}(y,s)}{(1+\alpha)^{1+a}}d\alpha,\\
\wt{J_{12}}(x,t)&\d:=\int_{R}^{t-x}\frac{d\alpha}{(1+\alpha)^{1+b}}\int_{\alpha}^{t-x}\frac{w^{-p}(-y,s)}{(1+\beta)^{1+a}}d\beta,\\
\wt{J_{13}}(x,t)&\d:=\int_{-R}^{R}\frac{d\beta}{(1+|\beta|)^{1+b}}\int_{|\beta|}^{t+x}\frac{w^{-p}(y,s)}{(1+\alpha)^{1+a}}d\alpha,\\
\wt{J_{14}}(x,t)&\d:=\int_{-R}^{R}\frac{d\alpha}{(1+|\alpha|)^{1+b}}\int_{|\alpha|}^{t-x}
\frac{w^{-p}(-y,s)}{(1+\beta)^{1+a}}d\beta,\\
\wt{J_{21}}(x,t)&\d:=\int_{-R}^{t-x}\frac{d\beta}{(1+|\beta|)^{1+b}}\int_{|\beta|}^{t+x}\frac{w^{-p}(y,s)}{(1+\alpha)^{1+a}}d\alpha,\\
\wt{J_{22}}(x,t)&\d:=\int_{-(t-x)_+}^{(t-x)_{+}}\frac{d\alpha}{(1+|\alpha|)^{1+b}}\int_{|\alpha|}^{t-x}
\frac{w^{-p}(-y,s)}{(1+\beta)^{1+a}}d\beta.\\
\end{array}
\]
\vskip10pt
\par\noindent
{\bf Case 1. $\v{a>0}$ and $\v{a+b>0}$, or $\v{a>0}$ and $\v{a+b=0}$, or $\v{a>0}$ and $\v{a+b<0}$. }

Because of $w(|x|,t)=1$, the desired estimates are established by the same manner as that of Case 1 of 
Lemma \ref{lem:apriori0}. We shall omit the detalis.
\vskip10pt
\par\noindent
{\bf Case 2. $\v{a=0}$ and $\v{b>0}$, or $\v{a=0}$ and $\v{b=0}$, or $\v{a=0}$ and $\v{b<0}$.} 

We shall estimate $\wt{J_{11}}$ in $D_{\rm Int}$. The definition of (\ref{weight}) yields
\begin{equation}
\label{J_11-est1}
\begin{array}{lll}
\wt{J_{11}}(x,t)&\le\d \int_{R}^{t-x}\frac{d\beta}{(1+\beta)^{1+b}}\int_{\beta}^{t+x}\frac{\log^p(\alpha+3R)d\alpha}{1+\alpha}\\
&\le\d C\log^p(t+x+3R)\int_{R}^{t-x}\frac{d\beta}{(1+\beta)^{1+b}}\int_{R}^{t+x}\frac{d\alpha}{1+\alpha}.
\end{array}
\end{equation}
Let  $a=0$ and $b>0$, or $a=0$ and $b=0$.
Making use of (\ref{est:alpha-beta}), we get 
\[
\int_{R}^{t-x}\frac{d\beta}{(1+\beta)^{1+b}}
\le
C\left\{
\begin{array}{lll}
1&\mbox{if} \ b>0,\\
\log(T+3R) &\mbox{if}\ b=0
\end{array}
\right.
\]
and 
\[
\int_{R}^{t+x}\frac{d\alpha}{1+\alpha}\le C\log(t+x+1)\le w^{-1}(x,t).
\]
It follows from (\ref{est:T1}) and (\ref{J_11-est1}) that
\[
\wt{J_{11}}(x,t)
\le
Cw^{-1}(x,t)\left\{
\begin{array}{lll}
\log^{p}(T+3R)&\mbox{if} \ b>0,\\
\log^{p+1}(T+3R) &\mbox{if}\ b=0.
\end{array}
\right.
\]

For the case of $a=0$ and $b<0$, we employ the same argument as in the proof of Lemma \ref{lem:apriori0}. By virtue of (\ref{est:alpha-beta}) and (\ref{est:T2}), we get  
\[
\begin{array}{lll}
\d\int_{R}^{t-x}\frac{d\beta}{(1+\beta)^{1+b}}\int_{R}^{t+x}\frac{d\alpha}{1+\alpha}
&\le&\d \int_{R}^{t-x}\frac{d\beta}{(1+\beta)^{1+b/2}}\int_{R}^{t+x}\frac{d\alpha}{(1+\alpha)^{1+b/2}}\\
&\le &\d C(T+R)^{-b}.
\end{array}
\]
Noticing that $\log^{p}(t+x+3R)\le \log^{p-1}(T+3R)w^{-1}(x,t)$ by (\ref{est:T1}), we obtain
\[
\wt{J_{11}}(x,t)\le C\log^{p-1}(T+3R)(T+R)^{-b}w^{-1}(x,t).
\]
Therefore, we get
\[
\wt{J_{11}}(x,t)\le Cw^{-1}(x,t)E_{2}(T)\quad \mbox{in}\ D_{\rm Int}.
\]

Next, we shall estimate $\wt{J_{12}}$. It follows from (\ref{weight}) that
\[
\wt{J_{12}}(x,t)\le \int_{R}^{t-x}\frac{d\alpha}{(1+\alpha)^{1+b}}\int_{\alpha}^{t+x}\frac{\log^p(\beta+3R)d\beta}{(1+\beta)^{1+a}}.
\]
The estimates of this integral are the same as those of $\wt{J_{11}}$, in which $\alpha$ and $\beta$ are replaced  with each other. Hence, we get
\[
\wt{J_{12}}(x,t)\le Cw^{-1}(x,t)E_{2}(T)\quad \mbox{in}\ D_{\rm Int}.
\]
Next, we shall estimate $\wt{J_{13}}$ and $\wt{J_{14}}$. The definition of (\ref{weight}) 
and (\ref{est:T1}) gives us 
\[
\begin{array}{ll}
\wt{J_{13}}(x,t)
&=\d\int_{-R}^{R}\frac{d\beta}{(1+|\beta|)^{1+b}}\int_{|\beta|}^{t+x}\frac{\log^{p}(\alpha+3R)d\alpha}{1+\alpha} \\
&\le \d C \log^p(T+3R)\int_{-R}^{t+x}\frac{d\alpha}{1+|\alpha|}
\end{array}
\]
and
\[
\begin{array}{ll}
\wt{J_{14}}(x,t)
&=\d \int_{-R}^{R}\frac{d\alpha}{(1+|\alpha|)^{1+b}}\int_{|\alpha|}^{t-x} 
\frac{\log^p(\beta+3R)d\beta}{1+\beta} \\
&\le \d C \log^p(T+3R)\int_{-R}^{t+x}\frac{d\alpha}{1+
|\beta|}.
\end{array}
\]
It follows from (\ref{ext-est1}), (\ref{TR-est}) and 
\[
\log^{p}(T+3R)\le C\log^{p-1}(T+3R)(T+R)^{-b}\quad \mbox{for}\ b<0
\]
that
\[
\wt{J_{13}}(x,t), \wt{J_{14}}(x,t)\le CE_{2}(T)w^{-1}(x,t).
\]

Next, we shall estimate $\wt{J_{21}}$ and $\wt{J_{22}}$ in $D_{\rm Ext} \cup D_{\rm Ori}$.
Since $|t-x|\le R$ holds for $(x,t)\in D_{\rm Ext}$, 
we have
\[
\begin{array}{ll}
\wt{J_{21}}(x,t)
&\le \d \int_{-R}^{R}\frac{d\beta}{(1+|\beta|)^{1+b}}\int_{-R}^{t+x}\frac{\log^{p}(\alpha+3R)d\alpha}{1+|\alpha|} \\
&\le \d C \log^p(T+3R)\int_{-R}^{t+x}\frac{d\alpha}{1+|\alpha|} 
\end{array}
\]
in $D_{\rm Ext}$, and 
\[
\begin{array}{ll}
\wt{J_{22}}(x,t)
&\le \d \int_{-R}^{R}\frac{d\beta}{(1+|\beta|)^{1+b}}\int_{-R}^{t-x}\frac{\log^{p}(\alpha+3R)d\alpha}{1+|\alpha|} \\
&\le \d C \log^p(T+3R)\int_{-R}^{t+x}\frac{d\alpha}{1+|\alpha|} 
\end{array}
\]
for $t-x>0$ in $D_{\rm Ext}$. Thus, the estimates of above integrals are the same as those of $\wt{J_{13}}$ and $\wt{J_{14}}$. 
Moreover, since $(x,t)\in D_{\rm Ori}$ is bounded, we have
\[
\wt{J_{21}}(x,t),\wt{J_{22}}(x,t) \le C \ \mbox{in $D_{\rm Ori}$}.
\]
Making use of $\log (t+|x|+3R)\ge 1$ and (\ref{TR-est}), we get the desired estimates. 
Hence, we obtain
\[
\wt{J_{21}}(x,t), \wt{J_{22}}(x,t)\le CE_{2}(T)w^{-1}(x,t) \quad \mbox{in}\ D_{\rm Ext} \cup D_{\rm Ori}.
\]
\vskip10pt
\par\noindent
{\bf Case 3. $\v{a<0}$ and $\v{b>0}$, or $\v{a<0}$ and $\v{b=0}$, or $\v{a<0}$ and $\v{b<0}$.} 

We shall show the estimate for $\wt{J_{11}}$ in $D_{\rm Int}$. The definition of $w$ in (\ref{weight}) and the trivial inequality (\ref{est:T2}) yield
\begin{equation}
\label{J_11-est2}
\begin{array}{lll}
\wt{J_{11}}(x,t)&\le&\d \int_{R}^{t-x}\frac{d\beta}{(1+\beta)^{1+b}}\int_{\beta}^{t+x}\frac{(\alpha+3R)^{-pa}d\alpha}{(1+\alpha)^{1+a}}\\
&\le &\d C(T+2R)^{-pa}\int_{R}^{t-x}\frac{d\beta}{(1+\beta)^{1+b}}\int_{R}^{t+x}\frac{d\alpha}{(1+\alpha)^{1+a}}.
\end{array}
\end{equation}
Making use of (\ref{est:alpha-beta}) for the above integral, we get
\[
\int_{R}^{t+x}\frac{d\alpha}{(1+\alpha)^{1+a}}\le w^{-1}(|x|,t)
\]
and
\[
\int_{R}^{t-x}\frac{d\beta}{(1+\beta)^{1+b}}
\le
C\left\{
\begin{array}{lll}
1&\mbox{if} \ b>0,\\
\log(T+3R) &\mbox{if}\ b=0,\\
(T+2R)^{-b}  &\mbox{if}\ b<0.
\end{array}
\right.
\]
Thus, (\ref{J_11-est2}) implies 
\[
\wt{J_{11}}(x,t)\le CE_{2}(T) w^{-1}(x,t)\quad \mbox{in}\ D_{\rm Int}.
\]

Next, we shall estimate $\wt{J_{12}}$. It follows from (\ref{weight}) that
\[
\wt{J_{12}}(x,t)\le \int_{R}^{t-x}\frac{d\alpha}{(1+\alpha)^{1+b}}\int_{\alpha}^{t+x}\frac{(\beta+3R)^{-pa}d\beta}{(1+\beta)^{1+a}}.
\]
The estimates of this integral are the same as those of $\wt{J_{11}}$, in which $\alpha$ and $\beta$ are replaced with each other. Hence, we get
\[
\wt{J_{12}}(x,t)\le Cw^{-1}(x,t)E_{2}(T)\quad \mbox{in}\ D_{\rm Int}.
\]
Next, we shall estimate $\wt{J_{13}}$ and $\wt{J_{14}}$. It follows from (\ref{weight}) 
and (\ref{est:T1}) that
\[\begin{array}{ll}
\wt{J_{13}}(x,t)
&=\d \int_{-R}^{R}\frac{d\beta}{(1+|\beta|)^{1+b}}\int_{|\beta|}^{t+x}\frac{(\alpha+3R)^{-pa}d\alpha}{(1+\alpha)^{1+a}} \\
&\le\d C (T+2R)^{-pa}\int_{-R}^{t+x}\frac{d\alpha}{
(1+|\alpha|)^{1+a}}
\end{array}
\]
and
\[
\begin{array}{ll}
\wt{J_{14}}(x,t)
&=\d\int_{-R}^{R}\frac{d\alpha}{(1+|\alpha|)^{1+b}}\int_{|\alpha|}^{t-x}
\frac{(\beta+3R)^{-pa}d\beta}{(1+\beta)^{1+a}} \\
&\le \d C (T+2R)^{-pa}\int_{-R}^{t+x}\frac{d\alpha}{(1+|\beta|)^{1+a}}.
\end{array}
\]
Hence, because of (\ref{ext-est1}) and (\ref{TR-est}), we get
\[
\wt{J_{13}}(x,t), \wt{J_{14}}(x,t)\le CE_{2}(T) w^{-1}(x,t)\quad \mbox{in}\ D_{\rm Int}.
\]

Next, we shall estimate $\wt{J_{21}}$ and $\wt{J_{22}}$ in $D_{\rm Ext} \cup D_{\rm Ori}$.
Since $|t-x|\le R$ holds for $(x,t)\in D_{\rm Ext}$,
we have
\[
\begin{array}{ll}
\wt{J_{21}}(x,t)
&\le \d \int_{-R}^{R}\frac{d\beta}{(1+|\beta|)^{1+b}}\int_{-R}^{t+x}
\frac{(\alpha+3R)^{-pa}d\alpha}{1+|\alpha|} \\
&\le\d C (T+2R)^{-pa}\int_{-R}^{t+x}\frac{d\alpha}
{(1+|\alpha|)^{1+a}} 
\end{array}
\]
in $D_{\rm Ext}$, and 
\[
\begin{array}{ll}
\wt{J_{22}}(x,t)
&\le\d \int_{-R}^{R}\frac{d\beta}{(1+|\beta|)^{1+b}}\int_{-R}^{t-x}\frac{(\alpha+3R)^{-pa}d\alpha}{1+|\alpha|}\\
&\le\d C (T+2R)^{-pa}\int_{-R}^{t+x}\frac{d\alpha}{(1+|\alpha|)^{1+a}} 
\end{array}
\]
for $t-x>0$ in $D_{\rm Ext}$. Thus, the estimates of above integrals are the same as those of $\wt{J_{13}}$ and $\wt{J_{14}}$. 
Moreover, since $(x,t)\in D_{\rm Ori}$ is bounded, we have
\[
\wt{J_{21}}(x,t),\wt{J_{22}}(x,t) \le C \ \mbox{in $D_{\rm Ori}$}.
\]
Making use of $t+|x|+3R\ge 1$ and (\ref{TR-est}), we get the desired estimates. 
Hence, we obtain
\[
\wt{J_{21}}(x,t), \wt{J_{22}}(x,t)\le CE_{2}(T)w^{-1}(x,t) \quad \mbox{in}\ D_{\rm Ext} \cup D_{\rm Ori}.
\]
Summing up all the estimates, we obtain (\ref{basic2}). Therefore the proof of (\ref{apriori2}) is now established.
\hfill
$\Box$

\vskip10pt
\par\noindent
{\bf Proofs of Theorem \ref{thm_global1} and Theorem \ref{thm_lower2}.}
\noindent
We shall employ the same argument as in the proof of Theorem 2.2 of \cite{KMT}.
 We consider an integral equation:
\begin{equation}
\label{integral2}
U=L(|U+\e u^0|^p)\quad \mbox{in}\ \R\times[0,T],
\end{equation}
where $L$ and $u^0$ are defined in (\ref{nonlinear}) and (\ref{linear}) respectively. 
Let us define a Banach space
\[
X_2:=\{U\in C(\R\times[0,T]): \ {\rm supp}\ U\subset \{|x|\le t+R\} \},
\]
which is equipped with a norm (\ref{norm2}). 
We shall construct a solution of the integral equation (\ref{integral2}) in a closed subspace 
\[
Y_2:=\{U\in X_2: \|U\|_2\le 2^{p+1}C_2\e^p\}
\]
of $X_2$, where $C_2$ is defined in (\ref{apriori1}).
Define a sequence of functions $\{U_n\}_{n\in \N}$ by 
\[
U_{n+1}=L(|U_n+\e u^0|^p),\quad U_1=0,
\]
Analogously to the proof of Theorem 1.2 in \cite{KMT} with $M:=C_2$, we can see that $U_n\in Y_2$ $(n\in\N)$ provided the inequality
\begin{equation}
\label{det:lifespan1}
2^{p^2+p}C_3C_2^{p-1}E_2(T)\e^{p(p-1)}\le 1
\end{equation}
holds, where $C_3$, $E_2(T)$ are defined in (\ref{apriori2}).
Moreover, $\{U_n\}$ is a Cauchy sequence in $Y_2$ provided the inequality 
\begin{equation}
\label{det:lifespan2}
3^{p-1}pC_32^{p^2+1}C_2^{p-1}E_2(T)\e^{p(p-1)}+3^{p-1}2pC_2D(T)\e^{p-1}\le 1
\end{equation}
holds, where $D(T)$ is the one in (\ref{apriori1}). 

When $a>0$, we can easily find $\e_1$ in Theorem \ref{thm_global1}, or $c$ and $\e_3$ in Theorem \ref{thm_lower2}, because of 
$D(T)=1$. 
When $a=0$ and $b>0$, the conditions (\ref{det:lifespan1}) and (\ref{det:lifespan2}) follow from 
\[
\log(T+3R)\le \wt{C_1}\e^{-(p-1)},
\]
where we set
\[
\begin{array}{ll}
\wt{C_1}:=\min&\!\!\!\! \{(2^{p^2+p}C_3C_2^{p-1})^{-1/p}, \\
& \ (2\cdot3^{p-1}pC_32^{p^2+1}C_2^{p-1})^{-1/p}, (2^23^{p-1}pC_2)^{-1}\}.
\end{array}
\]
When $a=b=0$, since $p(p-1)/(p+1)< p-1$, the conditions (\ref{det:lifespan1}) and (\ref{det:lifespan2}) follow from
\[
\log(T+3R)\le \wt{C_2}\e^{-p(p-1)/(p+1)},
\]
where we set
\[\begin{array}{ll}
\wt{C_2}:=\min&\!\!\!\!\{(2^{p^2+p}C_3C_2^{p-1})^{-1/(p+1)}, \\ &\ (2\cdot3^{p-1}pC_32^{p^2+1}C_2^{p-1})^{-1/(p+1)}, (2^23^{p-1}pC_2)^{-1}\}.
\end{array}
\]
When $a=0$ and $b<0$, the estimate in Theorem \ref{thm_lower2} is derived by 
\[
T\le \psi_{2}^{-1}(\wt{C_3}\e^{-p(p-1)}) \quad \ \mbox{for}\ 0<\e\le \e_3,
\]
where $\wt{C_3}:=(p2^{p^2+p-b+1}3^{p-1}C_3C_2^{p-1})^{-1}$. 
Here $\e_3$ has to satisfy
\[
R\le \psi_2^{-1}(\wt{C_3}\e_3^{-p(p-1)})\le \exp (3^{1-p}2^{-2}p^{-1}C_2^{-1}\e_3^{-(p-1)})-3R.
\]
When $a<0$ and $b>0$, 
the conditions (\ref{det:lifespan1}) and (\ref{det:lifespan2}) follow from
\[
T+2R\le \wt{C_4}\e^{-(p-1)/(-a)},
\]
where we set
\[
\begin{array}{ll}
\wt{C_4}:=\min&\!\!\!\!\{(2^{p^2+p}C_3C_2^{p-1})^{1/pa}, \\ & \ (2\cdot3^{p-1}pC_32^{p^2+1}C_2^{p-1})^{1/pa}, (2^23^{p-1}pC_2)^{1/a}\}.
\end{array}
\]
When $a<0$ and $b=0$, 
the estimate in Theorem \ref{thm_lower2} is derived by 
\[
T\le \psi_1^{-1}(\wt{C_5}\e^{-p(p-1)})\quad \ \mbox{for}\ 0<\e\le \e_3,
\]
where $\wt{C_5}:=(2^{p^2+3}3^{p-1-pa}pC_3C_2^{p-1})^{-1}$.
Here $\e_3$ has to satisfy
\[
R\le \psi_1^{-1}(\wt{C_5}\e_3^{-p(p-1)})\le (3^{1-p}2^{-2}p^{-1}C_2^{-1})^{-1/a}\e_3^{-(p-1)/(-a)}-2R.
\]
Finally, when $a<0$ and $b<0$, since $p(p-1)/(-pa-b)\le (p-1)/(-a)$, 
the conditions (\ref{det:lifespan1}) and (\ref{det:lifespan2}) follow from
\[
T+2R\le \wt{C_6}\e^{-p(p-1)/(-pa-b)},
\]
where we set
\[
\begin{array}{ll}
\wt{C_6}:=\min&\!\!\!\!\{(2^{p^2+p}C_3C_2^{p-1})^{1/(pa+b)}, \\ & \ (2\cdot3^{p-1}pC_32^{p^2+1}C_2^{p-1})^{1/(pa+b)}, (2^23^{p-1}pC_2)^{1/(pa+b)}\}.
\end{array}
\]
Therefore the proofs of Theorem \ref{thm_global1} with (\ref{initial_zero}) and Theorem \ref{thm_lower2} are completed.
\hfill
$\Box$
%kitamura_start%%%%%%%%%%%%%%%%%%%%%%%%%%%%%%%%%%%%%%%%%%%%%%%%%%%%%%%%%%%%%%%%%%%%%%

\section{Proofs of Theorem \ref{thm_upper1} and Theorem \ref{thm_upper2}}
\quad In this section, we shall investigate the upper bounds of the lifespan. We note that they are determined by point-wise estimates of the solution in the interior domain, $D_{\rm Int}$.
In fact, it follows from  (\ref{supp(f,g)}) and (\ref{linear}) that
\[
u(x,t)=\frac{\e}{2}\int_{\R}g(x)dx+L(|u|^p)(x,t)
\quad\mbox{for}\ (x,t)\in D_{\rm Int}.
\]
In this section, we assume that
\begin{equation}
\label{D}
(x,t)\in D:=D_{\rm Int}\cap\{x>0\}\cap\{t-x>R\}.
\end{equation}
Making use of Lemma \ref{weight_equivalent} and
introducing the characteristic coordinate by (\ref{coordinate}), we have that
\begin{equation}
\label{first}
u(x,t)\ge C_0
\int_R^{t-x}\frac{d\beta}{(1+\be)^{1+b}}\int_\beta^{t+x}\frac{|u(y,s)|^p}{(1+\al)^{1+a}}d\alpha
+J(x,t)
\end{equation}
and
\begin{equation}
\label{first2}
u(x,t)\ge C_0
\int_R^{t-x}\frac{d\alpha}{(1+\al)^{1+a}}\int_R^{\alpha}\frac{|u(y,s)|^p}{(1+\be)^{1+b}}d\beta
+J(x,t),
\end{equation}
where
\begin{equation}
\label{C_0}
C_0:=\frac{3^{-|a|}2^{-|b|}}{8\sqrt{2}}>0
\end{equation}
and
\begin{equation}
\label{J}
J(x,t):=
C_0\e^p\int_{-R}^R\frac{d\beta}{(1+|\be|)^{1+b}}\int_{|\be|}^{t+x}\frac{|u^0(y,s)|^p}{(1+\al)^{1+a}}d\alpha
+\frac{\e}{2}\int_{\R}g(x)dx.
\end{equation}
Employing this integral inequality, we shall estimate the lifespan from above.
We also use the following lemma. 
\begin{lem}\label{lemM_n}
For $\eta,\mu,C > 0$, Define $\{M_n\}$ by
\begin{equation}
\label{M_n}
M_{n+1}= C \eta^{-n} \mu^{-n} M_n^p .
\end{equation}
Then, 
\[
\log M_{n+1} \geq -\frac{1}{p-1}\log C+p^n\left\{\frac{1}{p-1}\log C-S_p\log(\eta\mu)+\log M_1\right\}
\]
holds, where 
\begin{equation}
\label{S_p}
S_p:=\sum_{j=1}^\infty\frac{j}{p^j}<\infty.
\end{equation}
\end{lem}
\par\noindent
{\bf Proof.}
\noindent
The definition of $M_n$ yields
\[
\log M_{n+1} = \log C - n \log (\eta\mu) + p \log M_n.
\]
Therefore, we obtain
\[
\begin{array}{ll}
\log M_{n+1}
&=(1+p+\cdots+p^{n-1})\log C\\
&\quad-\{n+p(n-1)+\cdots+p^{n-1}(n-n+1)\}\log (\eta\mu)+p^n\log M_1\\
&\d=\frac{p^n-1}{p-1}\log C-p^n\sum_{j=1}^n\frac{j}{p^j}\log (\eta\mu)+p^n\log M_1\\
&\d\ge-\frac{1}{p-1}\log C+p^n\left\{\frac{1}{p-1}\log C-S_p\log (\eta\mu)+\log M_1\right\}.
\end{array}
\]
The convergence of $S_p$ is trivial by d'Alembert's ratio test.
%It is clear from d'Alembert's ratio test that $S_p$ is finite.
\hfill$\Box$
%%% subsection 5.1 %%%

\subsection{Proof of Theorem \ref{thm_upper1}}
Let $u=u(x,t)\in C^2(\R\times[0,T))$ be a solution of (\ref{IVP}).
It follows from (\ref{positive_non-zero}), (\ref{first}) and (\ref{J}) that
\begin{equation}
\label{frame1}
u(x,t)\ge C_0
\int_R^{t-x}\frac{d\beta}{(1+\be)^{1+b}}\int_\beta^{t+x}\frac{|u(y,s)|^p}{(1+\al)^{1+a}}d\alpha+C_g\e
\end{equation}
and
\begin{equation}
\label{frame2}
u(x,t)\ge C_0
\int_R^{t-x}\frac{d\alpha}{(1+\al)^{1+a}}\int_R^{\alpha}\frac{|u(y,s)|^p}{(1+\be)^{1+b}}d\beta+C_g\e
\end{equation}
for $(x,t)\in D$, where
\[
C_g:=\frac{1}{2}\int_{\R}g(x)dx>0.
\]
\vskip10pt
\par\noindent
{\bf Case 1. $\v{a>0}$ and $\v{a+b<0}$, or $\v{a\leq0}$ and $\v{b<0}$.}

First, we consider only the case where $a>0$ and $a+b<0$.
Assume that an estimate
\begin{equation}
\label{case_a+b<0}
u(x,t)\ge M_n\left\{\frac{(t-x-R)^{2-b}}{(1+t-x)^{2+a}}\right\}^{a_n}
\quad\mbox{in}\ D
\end{equation}
holds, where $a_n\ge0$ and $M_n>0$.
The sequences $\{a_n\}$ and $\{M_n\}$ are defined later. 
Then it follows from (\ref{frame2}) and (\ref{case_a+b<0}) that
\[
u(x,t)\ge C_0M_n^p
\int_R^{t-x} \frac{d\al}{(1+\al)^{1+a}} \int_R^{\al} \frac{(\be-R)^{pa_n(2-b)} }{(1+\be)^{pa_n(2+a)+1+b}}d\be. 
\]
Because of 
\begin{equation}
\label{Cal1}
\beta-R \le 1+\beta
\end{equation}
and $b<0$, we get
\[
u(x,t)\ge \frac{C_0M_n^p}{(1+t-x)^{(pa_n+1)(2+a)}}
\int_R^{t-x} d\al \int_R^{\al}(\be-R)^{pa_n(2-b)-b}d\be
\]
which implies that
\begin{equation}
\label{uppera>0a+b<0}
u(x,t)\geq \frac{C_0M_n^p}{(2-b)^2(pa_n+1)^2}\left\{\frac{(t-x-R)^{2-b}}{(1+t-x)^{2+a}}\right\}^{pa_n+1}\quad\mbox{in}\ D.
\end{equation}
Therefore, if $\{a_n\}$ is defined by
\begin{equation}
\label{a_n}
a_{n+1}=pa_n+1,\ a_1=0,
\end{equation}
then (\ref{case_a+b<0}) holds for all $n\in\N$ as far as $M_n$ satisfies
\begin{equation}
%\label{M_n}
M_{n+1}\le\frac{C_0M_n^p}{(2-b)^2(pa_n+1)^2}.
\end{equation}
In view of (\ref{frame1}), we note that (\ref{case_a+b<0}) holds for $n=1$ with
\begin{equation}
\label{M_1}
M_1=C_g\e.
\end{equation}
Therefore, it follows from (\ref{a_n}) that
\begin{equation}
\label{a_n1}
a_n = \frac{p^{n-1}-1}{p-1}.
\end{equation}
According to $a_n<p^{n-1}/(p-1)$, (\ref{M_n}) and (\ref{M_1}), we define $\{M_n\}$ by
\[
M_{n+1} = C_1 p^{-2n} M_n^p,\quad M_1=C_g \e,
\]
where $C_1:=C_0(p-1)^2/(2-b)^2>0$. Hence, making use of Lemma \ref{lemM_n}, we reach to
\[
\log M_{n+1}\ge-\frac{1}{p-1}\log C_1+p^n\left\{\frac{1}{p-1}\log C_1-2S_p\log p+\log M_1\right\}.
\]
Therefore, we obtain
\[
u(x,t)\ge C_1^{-1/(p-1)}\left\{\frac{(t-x-R)^{2-b}}{(1+t-x)^{2+a}}\right\}^{-1/(p-1)}\exp\left(K_1(x,t)p^n\right)
\quad\mbox{in}\ D,
\]
where
\[
\begin{array}{ll}
K_1(x,t):=&\d \frac{1}{p-1}\log\left\{\frac{(t-x-R)^{2-b}}{(1+t-x)^{2+a}}\right\} \\
&\d +\frac{1}{p-1}\log C_1-2S_p\log p+\log M_1.
\end{array}
\]
If there exists a point $(x_0,t_0)\in D$ such that
\[
K_1(x_0,t_0)>0,
\]
we have
\[
u(x_0,t_0)=\infty
\]
by letting $n\rightarrow\infty$, so that $T<t_0$.
Let us set $2x_0=t_0$ and $6R < t_0$. Then $K_1(t_0/2,t_0)>0$ follows from
\begin{equation}
\label{con1}
2^{-2-a}(t_0/2-R)^{-a-b}C_1p^{-2S_p(p-1)}(C_g\e)^{p-1} > 1
\end{equation}
because the inequality 
\[
\frac{(t_0/2-R)^{2-b}}{(1+t_0/2)^{2+a}} \ge 2^{-2-a}(t_0/2-R)^{-a-b}
\]
holds for $6R<t_0$.
The condition (\ref{con1}) follows from
\[
%\label{blow-up_a+b<0}
t_0>4\{2^{-2-a}C_1p^{-2S_p(p-1)}C_g^{p-1}\}^{1/(a+b)}\e^{(p-1)/(a+b)}
\quad\mbox{for}\ 0<\e\le\e_4
\]
if we define $\e_4$ by
\[
2R=\{2^{-2-a}C_1p^{-2S_p(p-1)}(C_g\e_4)^{p-1}\}^{1/(a+b)}.
\]The proof for $a>0$ and $a+b<0$ is now completed.

Next, we consider the case where $a\leq0$ and $b<0$. Assume that
\[
%begin{equation}
%\label{case_a<0b<0}
u(x,t)\ge M_n\left\{\frac{(t-x-R)^{2-b-a}}{(1+t-x)^2}\right\}^{a_n}
\quad\mbox{in}\ D
\]
holds. Then, similarly to the above computations with (\ref{frame2}), we have
\[
u(x,t)\ge C_0M_n^p
\int_R^{t-x} \frac{d\al}{(1+\al)^{1+a}} \int_R^{\al} \frac{(\be-R)^{pa_n(2-b-a)} }{(1+\be)^{2pa_n+1+b}}d\be. 
\]
Making use of (\ref{Cal1}), we get 
\[
u(x,t)\geq \frac{C_0M_n^p}{(2-b-a)^2(pa_n+1)^2}\left\{\frac{(t-x-R)^{2-b-a}}{(1+t-x)^2}\right\}^{pa_n+1}\!\!\!
\quad\mbox{in}\ D.
\]
This estimate is almost same as (\ref{uppera>0a+b<0}). Therefore, the same lifespan estimate is obtained except for the included constant independent of $\e$.
\hfill$\Box$
\vskip10pt
\par\noindent
{\bf Case 2. $\v{a<0}$ and $\v{b=0}$.} 

Assume that an estimate
\begin{equation}
\label{case_a<0b=0}
u(x,t) \ge M_n\left\{ \frac{(2x)^{1-a}}{1+t+x}\log\left(\frac{1+t-x}{1+R}\right) \right\}^{a_n}
\quad\mbox{in}\ D
\end{equation}
holds, where $a_n$ is the one in (\ref{a_n1}).
In this case, $M_n$ is defined by (\ref{M_n}) with $C=C_2:=C_0(p-1)^2(1-a)^{-1}>0$ and $\eta=\mu=p$. 
We note that the same argument as in Case 1 can be applicable to also this case.
It follows from (\ref{frame1}) and (\ref{case_a<0b=0}) that
\[
u(x,t)\ge C_0M_n^p
\int_R^{t-x}\log^{pa_n}\!\left(\frac{1+\be}{1+R}\right)\frac{d\be}{1+\beta}\int_\be^{t+x} \frac{(\al-\be)^{pa_n(1-a)} }{(1+\al)^{1+a+pa_n}}d\al.
\]
For the $\al$-integral, we obtain
\begin{equation}
\label{al-est1}
\begin{array}{ll}
&\d \int_\be^{t+x}\frac{(\al-\be)^{pa_n(1-a)}d\al}{(1+\al)^{1+a+pa_n}}\\
&\d \ge \frac{1}{(1+t+x)^{pa_n+1}}\int_\be^{t+x}(\al-\be)^{pa_n(1-a)-a}d\al\\
&\d = \frac{(t+x-\be)^{(pa_n+1)(1-a)}}{(pa_n+1)(1-a)(1+t+x)^{pa_n+1}}.
\end{array}
\end{equation}
Replacing $\beta$ by $(t-x)$ in the quantity above, we get 
\[
u(x,t)\geq \frac{C_0M_n^p}{(1-a)(pa_n+1)^2}\left\{\frac{(2x)^{1-a}}{1+t+x}\log\left(\frac{1+t-x}{1+R}\right) \right\}^{pa_n+1}\ \mbox{in}\ D.
\]
Therefore, we obtain 
\[
u(x,t)\ge C_2^{-1/(p-1)}\left\{\frac{(2x)^{1-a}}{1+t+x}\log\left(\frac{1+t-x}{1+R}\right) \right\}^{-1/(p-1)}\exp\left(K_2(x,t)p^n\right)
\]
in $D$, where $K_2$ is defined by 
\[
\begin{array}{ll}
K_2(x,t):=
&\d \frac{1}{p-1}\log\left\{\frac{(2x)^{1-a}}{1+t+x}\log\left(\frac{1+t-x}{1+R}\right) \right\}\\
&\d +\frac{1}{p-1}\log C_2-2S_p\log p+\log M_1
\end{array}
\]
and $S_p$ is the one in (\ref{S_p}). %lemma \ref{lemM_n}. 
%When $2x=t$ and $2+t \ge 16R^2$, 
%\[
%\frac{(2x)^{1-a}}{1+t+x} = \frac{t^{1-a}}{1+3t/2} \ge \frac{t^{-a}}{3}
%\]
%and
%\[
%\log\left(\frac{1+t-x}{1+R}\right) \ge \log\left(\frac{1+t/2}{2R}\right) \ge \frac{1}{2}\log(2+t)
%\]
%hold. 
If there exists $(x_0,t_0)\in D$ such that
%making use the same argument as in {\bf Case1}, 
$K_2(x_0,t_0) > 0$, then we have $u(x_0,t_0) = \infty$ as before.
Set $2x_0=t_0$ and $2+t_0\ge16R^2$. Then $K_2(x_0,t_0) > 0$ follows from
\[
\phi(t_0)=t_0^{-a}\log(2+t_0)>3 \cdot 2 \{C_2p^{-2S_p(p-1)}C_g^{p-1}\}^{-1}\e^{-(p-1)}
\]
for $0<\e\le \e_4$, where $\e_4$ is defined by 
\[
\phi^{-1}(3 \cdot 2 \{C_2p^{-2S_p(p-1)}C_g^{p-1}\}^{-1}\e_4^{-(p-1)})= 16R^2-2
\]
because of 
\[
\frac{(2x_0)^{1-a}}{1+t_0+x_0} = \frac{t_0^{1-a}}{1+3t_0/2} \ge \frac{t_0^{-a}}{3}
\]
and
\[
\log\left(\frac{1+t_0-x_0}{1+R}\right) \ge \log\left(\frac{1+t_0/2}{2R}\right) \ge \frac{1}{2}\log(2+t_0).
\]
The proof for $a<0$ and $b=0$ is now completed.
\hfill$\Box$
\vskip10pt
\par\noindent
{\bf Case 3. $\v{a<0}$ and $\v{b>0}$.} 

Assume that an estimate
\begin{equation}
\label{case_a<0b>0}
u(x,t) \ge M_n\left\{ \frac{(2x)^{1-a}}{1+t+x}\frac{t-x-R}{(1+t-x)^{1+b}} \right\}^{a_n}
\quad\mbox{in}\ D
\end{equation}
holds, where $a_n$ is the one in (\ref{a_n1}).
In this case, $M_n$ is defined by (\ref{M_n}) with $C=C_3:=C_0(p-1)^2(1-a)^{-1}>0$ and $\eta=\mu=p$.

Making use of (\ref{frame1}) and (\ref{case_a<0b>0}), we have
\[
u(x,t)\ge C_0M_n^p
\int_R^{t-x} \frac{(\be-R)^{pa_n}d\be}{(1+\be)^{(pa_n+1)(1+b)}}
\int_\be^{t+x}\frac{(\al-\be)^{pa_n(1-a)}d\al}{(1+\al)^{1+a+pa_n}}.
\]
It follows from 
\[
(pa_n+1)(1+b)>0, \ pa_n+1>0 \ \mbox{and} \ 1+\al\ge (\al-\be)
\]
that
\[
\begin{array}{lll}
u(x,t)&\ge\d \frac{C_0M_n^p}{(1+t-x)^{(pa_n+1)(1+b)}(1+t+x)^{pa_n+1}}\\
 &\d\quad\times\int_R^{t-x} (\be-R)^{pa_n}d\be
\int_\be^{t+x}(\al-\be)^{pa_n(1-a)-a}d\al.
\end{array}
\]
Hence (\ref{al-est1}) implies that
\[
u(x,t)\geq \frac{C_0M_n^p}{(1-a)(pa_n+1)^2}\left\{\frac{(2x)^{1-a}}{1+t+x} \frac{t-x-R}{(1+t-x)^{1+b}} \right\}^{pa_n+1}\quad\mbox{in}\ D.
\]
Therefore, we obtain 
\[
u(x,t)\ge C_3^{-1/(p-1)}\left\{\frac{(2x)^{1-a}}{1+t+x}\frac{t-x-R}{(1+t-x)^{1+b}} \right\}^{-1/(p-1)}\exp\left(K_3(x,t)p^n\right)
\]
in $D$, where $K_3$ is defined by 
\[
\begin{array}{ll}
K_3(x,t):=
&\d \frac{1}{p-1}\log\left\{ \frac{(2x)^{1-a}}{1+t+x} \frac{t-x-R}{(1+t-x)^{1+b}}\right\}\\
&\d +\frac{1}{p-1}\log C_3-2S_p\log p+\log M_1
\end{array}
\]
and $S_p$ is the one in (\ref{S_p}). Let us fix $(x_0,t_0) \in D$ such that 
\[
t_0 - x_0 = R+1\ \mbox{and} \ t_0 > 2R.
\] 
Then we have 
\[
3R(t_0-R-1) \geq (R-1)(2t_0-R).
\]
Hence $K_3(x_0,t_0)>0$ follows from 
%\begin{equation}
%\label{blow-up_a<0b>0}
\[
(t_0-R-1)^{-a}>\frac{2^{a-1}\cdot 3R (2+R)^{1+b}}{(R-1)C_3p^{-2S_p(p-1)}C_g^{p-1}}\e^{-(p-1)}.
\]
This inequality provides us the desired estimate as before.
The proof for $a<0$ and $b>0$ is now completed.
\hfill$\Box$
\vskip10pt
\par\noindent
{\bf Case 4. $\v{a=0}$ and $\v{b>0}$.} 

Assume that an estimate
\begin{equation}
\label{case_a=0b>0}
u(x,t) \ge M_n\left\{\log\left(\frac{1+t+x}{1+t-x}\right)\frac{t-x-R}{(1+t-x)^{1+b}} \right\}^{a_n}
\quad\mbox{in}\ D
\end{equation}
holds, where $a_n$ is the one in (\ref{a_n1}). In this case
$M_n$ is defined by (\ref{M_n}) with
$C=C_4:=C_0(p-1)^2>0$ and $\eta=\mu=p$.

Making use of (\ref{frame1}) and (\ref{case_a=0b>0}), we have
\[
\begin{array}{llll}
u(x,t)&\ge\d C_0M_n^p
\int_R^{t-x}\frac{(\be-R)^{pa_n}}{(1+\be)^{(pa_n+1)(1+b)}}d\be 
\int_\be^{t+x} \log^{pa_n}\left(\frac{1+\al}{1+\be}\right)\frac{d\al }{(1+\al)}\\
&\ge\d \frac{C_0M_n^p}{(1+t-x)^{(pa_n+1)(1+b)}}  \\
&\quad\times\d \int_R^{t-x}(\be-R)^{pa_n}d\be 
\int_\be^{t+x} \log^{pa_n}\left(\frac{1+\al}{1+\be}\right)\frac{d\al }{(1+\al)}.
\end{array}
\]
For the $\al$-integral, we have that
\[
\int_\be^{t+x}\log^{pa_n}\left(\frac{1+\al}{1+\be}\right)\frac{d\al}{1+\al} = \frac{1}{pa_n+1} \log^{pa_n+1}\left(\frac{1+t+x}{1+\be}\right),
\]
which gives us
\[
\begin{array}{ll}
u(x,t)&\d \ge \frac{C_0M_n^p}{(pa_n+1)(1+t-x)^{(pa_n+1)(1+b)}}\\
&\quad\times \d \log^{pa_n+1}\left(\frac{1+t+x}{1+t-x}\right)
\int_R^{t-x}(\be-R)^{pa_n}d\be,
\end{array}
\]
which implies that
\[
u(x,t)\geq \frac{C_0M_n^p}{(pa_n+1)^2}\left\{\log\left(\frac{1+t+x}{1+t-x}\right)\frac{t-x-R}{(1+t-x)^{1+b}}\right\}^{pa_n+1}\quad \mbox{in}\ D.
\]
Therefore, we have 
\[
u(x,t)\ge C_4^{-1/(p-1)}\left\{\log\left(\frac{1+t+x}{1+t-x}\right)\frac{t-x-R}{(1+t-x)^{1+b}} \right\}^{-1/(p-1)}\!\!\!\!\exp\left(K_4(x,t)p^n\right)
\]
in $D$, where $K_4$ is defined by 
\[
\begin{array}{ll}
K_4(x,t):=
&\d \frac{1}{p-1}\log\left\{\log\left(\frac{1+t+x}{1+t-x}\right)\frac{t-x-R}{(1+t-x)^{1+b}}\right\}\\
&\d +\frac{1}{p-1}\log C_4-2S_p\log p+\log M_1
\end{array}
\]
and $S_p$ is the one in (\ref{S_p}). As in Case 3, let us fix $(x_0,t_0) \in D$ such that 
\[
t_0 - x_0 = R+1\ \mbox{and} \ t_0 > (R+2)^2.
\] 
Then we have
\[
\log\left(\frac{2t_0-R}{R+2}\right)\ge \log\left(\frac{t_0}{R+2}\right) \ge \frac{1}{2}\log t_0.
\]
Hence, $K_4(x_0,t_0)>0$ follows from 
%\begin{equation}
%\label{blow-up_a=0b>0}
\[
\log t_0>2(1+R)^{1+b}\{C_4p^{-2S_p(p-1)}C_g^{p-1}\}^{-1}\e^{-(p-1)}.
\]
This inequality provides us with the desired estimate as before.
The proof for $a<0$ and $b>0$ is now completed.
\hfill$\Box$
\vskip10pt
\par\noindent
{\bf Case 5. $\v{a=0}$ and $\v{b=0}$.} 

Assume that an estimate
\begin{equation}
\label{case_a=0b=0}
u(x,t) \ge M_n\left\{\log^2\left(\frac{1+t-x}{1+R}\right)\right\}^{a_n}
\quad\mbox{in}\ D
\end{equation}
holds, where $a_n$ is the one in (\ref{a_n1}). In this case 
$M_n$ is defined by (\ref{M_n}) with $C=C_5:=C_0(p-1)^2\cdot 4^{-1}>0$ and $\eta=\mu=p$.
Making use of (\ref{frame2}) and (\ref{case_a=0b=0}), we have
\[
\begin{array}{llll}
u(x,t)
&\ge&\d C_0M_n^p
\int_R^{t-x}\frac{d\al}{1+\al}\int_R^\al \log^{2pa_n}\left(\frac{1+\be}{1+R}\right)\frac{d\be}{1+\be}\\
&=& \d  \frac{C_0M_n^p}{2pa_n+1}
\int_R^{t-x}\log^{2pa_n+1}\left(\frac{1+\al}{1+R}\right)\frac{d\al}{1+\al}\\
&\ge&\d \frac{C_0M_n^p}{4(pa_n+1)^2}\log^{2(pa_n+1)}\left(\frac{1+t-x}{1+R}\right)
\end{array}
\]
%The above inequation can be obtained by cutting the integral domain into triangles and exchanging the integral order. 
in $D$.
Therefore, we have
\[
u(x,t)\ge C_5^{-1/(p-1)}\left\{\log^2\left(\frac{1+t-x}{1+R}\right) \right\}^{-1/(p-1)}\exp\left(K_5(x,t)p^n\right)
\]
in $D$, where $K_5$ is defined by
\[
\begin{array}{ll}
K_5(x,t):=
&\d \frac{1}{p-1}\log\left\{\log^2\left(\frac{1+t-x}{1+R}\right)\right\}\\
&\d +\frac{1}{p-1}\log C_5-2S_p\log p+\log M_1
\end{array}
\]
and $S_p$ is the one in (\ref{S_p}). As in Case 1, let us fix $(x_0,t_0) \in D$ such that 
\[
t_0=2x_0 \mbox{ and } t_0 > 4(1+R)^2.
\]
Then we have
\[
\log\frac{1+t_0/2}{1+R} \ge \log\frac{t_0}{2(1+R)}  \ge \frac{1}{2}\log t_0,
\]
Hence $K_5(x_0,t_0)>0$ follows from 
\[
\log t_0>2\{C_5p^{-2S_p(p-1)}C_g^{p-1}\}^{-1/2}\e^{-(p-1)/2}
\]
which is the desired estimate as before.
The proof for $a=0$ and $b=0$ is now completed.
\hfill$\Box$
\vskip10pt
\par\noindent
{\bf Case 6. $\v{a>0}$ and $\v{a+b=0}$.} 

In this case, we employ so-called \lq\lq slicing method'' which was introduced by Agemi, Kurokawa and Takamura \cite{AKT00}. 
Let us set
\begin{equation}
\label{D_l_n}
(x,t)\in D_{n}:=D_{\rm Int}\cap\{x>0\}\cap\{t-x>l_nR\}, \quad l_n:=\sum_{i=0}^n\left(\frac{1}{2}\right)^i.
\end{equation}
We note that $D_{n+1} \subset D_n$ for all $n \in \N$.  
Assume that an estimate
\begin{equation}
\label{case_a>0a+b=0}
u(x,t) \ge M_n\left\{\log\left(\frac{1+t-x}{1+l_nR}\right)\right\}^{a_n}
\quad\mbox{in}\ D_{n}
\end{equation}
holds, where $a_n$ is the one in (\ref{a_n1}).  The sequence $\{M_n\}$ with $M_n>0$ is defined later. 
Then (\ref{frame1}) and (\ref{case_a>0a+b=0}) imply that
\[
u(x,t)\ge C_0M_n^p
\int_{l_nR}^{t-x}\frac{d\al}{(1+\al)^{1+a}}\int_{l_nR}^\al \log^{pa_n}\left(\frac{1+\be}{1+l_nR}\right)\frac{d\be}{(1+\be)^{1+b}}.
\]
Let $(x,t)\in D_{n+1}$. Then, we get
\[
u(x,t)\ge C_0M_n^p
\int_{l_{n+1}R}^{t-x}\frac{d\al}{(1+\al)^{1+a}}\int_{l_nR}^\al\log^{pa_n}\left(\frac{1+\be}{1+l_nR}\right) \frac{d\be}{(1+\be)^{1+b}}.
\]
Note that
\[
l_nR <\lambda_n(\al)<\al \ \mbox{ for all }\ n \in \N
\]
holds, where
\[
\lambda_n(\al):=\frac{(1+\al)(1+l_nR)}{1+l_{n+1}R}-1.
\]
Then we have that
\[
\begin{array}{llll}
u(x,t)&\ge\d C_0M_n^p
\int_{l_{n+1}R}^{t-x}\frac{d\al}{(1+\al)^{1+a}}\int_{\lambda_n}^\al\log^{pa_n}\left(\frac{1+\be}{1+l_nR}\right) \frac{d\be}{(1+\be)^{1+b}}\\
&\ge\d C_0M_n^p
\int_{l_{n+1}R}^{t-x}\log^{pa_n}\left(\frac{1+\al}{1+l_{n+1}R}\right)\frac{d\al}{(1+\al)^{1+a}}\int_{\lambda_n}^\al \frac{d\be}{(1+\be)^{1+b}}.
\end{array}
\]
Due to $1<l_n<2$ for all $n \in \N$ and $R>1$, the $\beta$-integral is estimated as follows: 
\[
\begin{array}{llll}
\d \int_{\lambda_n}^\al \frac{d\be}{(1+\be)^{1+b}}&\ge\d \left(\frac{1+l_{n+1}R}{1+l_{n}R}\right)^{b}\frac{1}{(1+\al)^{b+1}}\int_{\lambda_n}^\al d\be\\
&=\d  \frac{(1+l_{n+1}R)^{b-1}}{(1+l_{n}R)^b}\cdot \frac{(l_{n+1}-l_n)R}{(1+\al)^b}\\
&\ge\d \frac{3^{b-1}\cdot 2^{-(n+1)}}{(1+\al)^b}.
\end{array}
\]
Hence, we get
\[
\begin{array}{llll}
u(x,t)
&\d \ge C_0M_n^p3^{b-1}\cdot 2^{-(n+1)}\int_{l_{n+1}R}^{t-x}\log^{pa_n}\left(\frac{1+\al}{1+l_{n+1}R}\right)\frac{d\al}{(1+\al)^{1+a+b}} \\
&\d \ge\frac{C_0M_n^p3^{b-1}\cdot 2^{-(n+1)}}{pa_n+1}\log^{pa_n+1}\left(\frac{1+t-x}{1+l_{n+1}R}\right) \quad \mbox{in}\ D_{n+1}.
\end{array}
\] 
Therefore, we can employ the same argument as in Case 1 with $a_n$ in (\ref{a_n1})
and $M_n$ defined by
\[
M_{n+1}=C_6 2^{-n}p^{-n}M_n^p , \quad M_1 = C_g \e,
\]
where $C_6:=C_0(p-1)/2\cdot 3^{1-b}>0$.
Here, we assume $\d t-x> 2R = \lim_{n \to \infty} l_n R$. Making use of Lemma \ref{lemM_n} with $C=C_6$, $\eta = 2$ and $\mu = p$, we obtain
\[
u(x,t)\ge C_6^{-1/(p-1)}\left\{\log\left(\frac{1+t-x}{1+2R}\right)\right\}^{-1/(p-1)}\exp\left(K_6(x,t)p^n\right)
\quad\mbox{in}\ D_\infty,
\]
where
\[
\begin{array}{ll}
K_6(x,t):=&
 \d \frac{1}{p-1}\log\left\{\log\left(\frac{1+t-x}{1+2R}\right)\right\} \\
 &\d+\frac{1}{p-1}\log C_6-S_p\log (2p)+\log M_1
 \end{array}
\]
and $S_p$ is the one in (\ref{S_p}). As in Case 1, let us fix $(x_0,t_0) \in D_\infty$ such that 
\[
t_0=2x_0\ \mbox{and} \ t_0> 4(1+2R)^2.
\]
Hence $K_6(x_0,t_0)>0$ follows from 
\[
\log t_0>2\{C_6(2p)^{-S_p(p-1)}C_g^{p-1}\}^{-1}\e^{-(p-1)}
\]
which is the desired estimate as before.
The proof for $a>0$ and $a+b=0$ is now completed. 
Therefore the proof of Theorem \ref{thm_upper1} finishes.
\hfill$\Box$

%%% subsection5.2 %%%

\subsection{Proof of Theorem \ref{thm_upper2}}
The proof is almost similar to the one of Theorem \ref{thm_upper1}. 
Let $u=u(x,t)\in C^2(\R\times[0,T))$ be a solution of (\ref{IVP}). Since the assumption on the initial data in (\ref{positive_zero}) yields 
\begin{equation}
\label{lower-bound_linear}
u^0(x,t)=\frac{1}{2}\{f(x+t)+f(x-t)\}\ge\frac{1}{2}f(x-t)
\quad \mbox{for} \ (x,t)\in\R\times[0,\infty).
\end{equation}
It follows from (\ref{first}), (\ref{first2}), (\ref{J}) and
\[
u(x,t)\ge\e u^0(x,t)\ge\frac{\e}{2}f(x-t) \quad \mbox{for} \ (x,t)\in\R\times[0,T)
\]
that
\begin{equation}
\label{frame3}
u(x,t)\ge C_0
\int_R^{t-x}\frac{d\beta}{(1+\be)^{1+b}}\int_\beta^{t+x}\frac{|u(y,s)|^p}{(1+\al)^{1+a}}d\alpha
+\frac{C_0}{2^p}\e^pJ'(x,t)
\end{equation}
and
\begin{equation}
\label{frame4}
u(x,t)\ge C_0
\int_R^{t-x}\frac{d\alpha}{(1+\al)^{1+a}}\int_R^{\alpha}\frac{|u(y,s)|^p}{(1+\be)^{1+b}}d\beta
+\frac{C_0}{2^p}\e^pJ'(x,t)
\end{equation}
for $(x,t)\in D$, where $D$ and $C_0$ are defined in (\ref{D}) and (\ref{C_0}) respectively,  and
\begin{equation}
\label{J'}
J'(x,t):=\int_{-R}^R\frac{f(-\beta)^p}{(1+|\be|)^{1+b}}d\be\int_{|\be|}^{t+x}\frac{d\al}{(1+\al)^{1+a}}.
\end{equation}
For (\ref{J'}), we may assume that
\begin{equation}
\label{positive_zero_add}
f(z) \not\equiv 0\quad \mbox{in}\quad z \in (-R,0)
\end{equation} 
without loss of generality. Because, if not, we have to assume $f(z) \not\equiv 0$ in $z \in (0,R)$ and to change the definition of $D$ in which $x>0$ is replaced with $x<0$. For such a case, we obtain all the estimates below with $-x$ instead of $x$.
In fact, taking $f(x+t)$ instead of $f(x-t)$ in (\ref{lower-bound_linear}),
we have, instead of (\ref{frame2}), that
\[
u(x,t)\ge C_0
\int_R^{t+x}\frac{d\alpha}{(1+\al)^{1+b}}\int_\alpha^{t-x}\frac{|u(y,s)|^p}{(1+\be)^{1+a}}d\beta
+\frac{C_0}{2^p}\e^pJ''(x,t),
\]
where
\[
J''(x,t):=\int_{-R}^R\frac{f(\alpha)^p}{(1+|\al|)^{1+b}}d\alpha\int_{|\al|}^{t-x}\frac{d\be}{(1+\be)^{1+a}}.
\]

\vskip10pt
\par\noindent
{\bf Case 1. $\v{a>0}$ and $\v{a+b<0}$.}

Since $J'$ in (\ref{J'}) is estimated from below by
\[
\begin{array}{ll}
\d \int_0^R \frac{f(-\be)^p}{(1+\be)^{1+b}} d\be \int_\be^{t-x} \frac{d\al}{(1+\al)^{1+a}}
&\d \ge \int_0^R \frac{f(-\be)^p}{(1+\be)^{1+b}} d\be \int_\beta^{R} \frac{d\al}{(1+\al)^{1+a}}\\
&\d \ge \frac{1}{(1+R)^{1+a}} \int_0^R\frac{(R-\beta)f(-\be)^p}{(1+\beta)^{1+b}}d\be,  \\
\end{array}
\]
it follows from (\ref{frame4}) and (\ref{positive_zero_add}) that
\[
u(x,t)\ge C_0
\int_{R}^{t-x}d\al\int_{R}^{\al} \frac{|u(y,s)|^p}{(1+\al)^{1+a}(1+\be)^{1+b}}d\be
+C_{f,1}\e^p \quad \mbox{in}\ D,
\]
where
\begin{equation}
\label{Cf1}
C_{f,1}:=\frac{C_0}{2^p(1+R)^{1+a}}\int_0^R\frac{(R-\beta)f(-\be)^p}{(1+\beta)^{1+b}}d\be>0.
\end{equation}
Therefore, we can employ the same argument as Case 1 of Theorem \ref{thm_upper1} in which the constant $C_g\e$ in (\ref{frame2}) is simply replaced with $C_{f,1}\e^p$.
\hfill$\Box$
\vskip10pt
\par\noindent
{\bf Case 2. $\v{a=0}$ and $\v{b<0}$.}

Since $J'$ in (\ref{J'}) is estimated from below by
\[
\int_0^R \frac{f(-\be)^p}{(1+\be)^{1+b}} d\be \int_\be^{t-x} \frac{d\al}{1+\al}
 \ge \log\left(\frac{1+t-x}{1+R}\right)\int_0^R \frac{f(-\be)^p}{(1+\be)^{1+b}} d\be,
\]
it follows from (\ref{frame3}) and (\ref{positive_zero_add}) that
\begin{equation}
\label{frame2_2}
\begin{array}{ll}
u(x,t)&\d \ge C_0
\int_{R}^{t-x}\frac{d\be}{(1+\be)^{1+b}}\int_\be^{t-x} \frac{|u(y,s)|^p}{1+\al}d\al \\
& \d +C_{f,2}\e^p\log\left(\frac{1+t-x}{1+R}\right)
\end{array}
\end{equation}
for $(x,t)\in D$, where
\begin{equation}
\label{Cf2}
C_{f,2}:=\frac{C_0}{2^p}\int_0^R\frac{f(-\be)^p}{(1+\be)^{1+b}}d\beta>0. 
\end{equation}
We employ the \lq\lq slicing method'' again. Assume that an estimate
\begin{equation}
\label{case_a=0b<0g=0}
u(x,t)\ge M_n\left(\frac{1+t-x}{1+l_nR}\right)^{-b\cdot a_n}\left\{\log\left(\frac{1+t-x}{1+l_nR}\right)\right\}^{b_n}
\quad\mbox{in}\ D_n
\end{equation}
holds, where $a_n \ge 0$, $b_n > 0$ and  $M_n>0$. Here, $D_n$ and $l_n$ are defined in (\ref{D_l_n}).
The sequences $\{a_n\}$, $\{b_n\}$ and $\{M_n\}$ are defined later. 
Then it follows from (\ref{frame2_2}) and (\ref{case_a=0b<0g=0}) that
\[
\begin{array}{ll}
u(x,t) \\
 \d \ge C_0M_n^p
\int_{l_nR}^{t-x}\left(\frac{1+\be}{1+l_nR}\right)^{-b\cdot pa_n}\log^{pb_n}\left(\frac{1+\be}{1+l_nR}\right)\frac{d\be}{(1+\be)^{1+b}} \int_{\be}^{t-x}\frac{d\al}{1+\al} \\
 \d \ge \frac{C_0M_n^p}{(1+t-x)^2}
\int_{l_nR}^{t-x} \frac{t-x-\be}{(1+\be)^b}\left(\frac{1+\be}{1+l_nR}\right)^{-b\cdot pa_n}\log^{pb_n}\left(\frac{1+\be}{1+l_nR}\right)d\be.
\end{array}
\]
Let $(x,t)\in D_{n+1}$. Note that
\[
l_nR <\lambda_n<t-x \quad \mbox{for all} \ n \in \N
\]
holds, where
\[
\lambda_n:=\lambda_n(t-x)=\frac{(1+t-x)(1+l_nR)}{1+l_{n+1}R}-1.
\]
Then we have
\[
\begin{array}{lllll}
u(x,t)
&\d \ge \frac{C_0M_n^p}{(1+t-x)^2} \\
&\d \quad \times \int_{\lambda_n}^{t-x} \frac{t-x-\be}{(1+\be)^b}\left(\frac{1+\be}{1+l_nR}\right)^{-b\cdot pa_n}\log^{pb_n}\left(\frac{1+\be}{1+l_nR}\right)d\be \\
&\d \ge \frac{C_0M_n^p}{(1+t-x)^2}
\left(\frac{1+l_nR}{1+l_{n+1}R}(1+t-x)\right)^{-b}\left(\frac{1+t-x}{1+l_{n+1}R}\right)^{-b\cdot pa_n} \\ &\quad \d \times\log^{pb_n}\left(\frac{1+t-x}{1+l_{n+1}R}\right)\int_{\lambda_n}^{t-x} (t-x-\be)d\be. \\
\end{array}
\]
It follows from
\[
\int_{\lambda_n}^{t-x} (t-x-\be)d\be=\frac{1}{2}\left\{\frac{(l_{n+1}-l_n)R}{1+l_{n+1}R}\right\}^2(t-x+1)^2\ge \frac{2^{-2n}}{72}(1+t-x)^2
\]
that
\[
u(x,t)\ge \frac{C_02^{-2n}M_n^p}{72}
\left(\frac{1+t-x}{1+l_{n+1}R}\right)^{-b( pa_n+1)}\log^{pb_n}\left(\frac{1+t-x}{1+l_{n+1}R}\right) \quad \mbox{in}\ D_{n+1}.
\]
Hence, (\ref{case_a=0b<0g=0}) holds for all $n\in\N$ provided
\[
\left\{
\begin{array}{lc}
a_{n+1}=pa_n+1,\ & a_1=0\\
b_{n+1}=pb_n,\ & b_1=1
\end{array}
\right.
\]
and 
\[
M_{n+1} \le C_7 2^{-2n}M_n^p, \ M_1=C_{f,2}\e^p,
\] 
where $C_7 := C_0/72>0$.
Therefore, we define $a_n$ is the one in (\ref{a_n1}), $b_n = p^{n-1}$ and $M_n$ as above. Making use of Lemma \ref{lemM_n} with $C=C_7$ and $\eta = \mu =2$, for $(x,t)\in D_\infty$, we obtain
\[
u(x,t)\ge C_7^{-1/(p-1)}(1+t-x)^{b/(p-1)}\exp\left(K_7(x,t)p^{n-1}\right),
\]
where
\[
\begin{array}{ll}
 K_7(x,t):=&\d \frac{1}{p-1}\log\left\{\left(\frac{1+t-x}{1+2R}\right)^{-b}\log^{p-1}\left(\frac{1+t-x}{1+2R}\right) \right\} \\
&\d +\frac{1}{p-1}\log C_7-2S_p\log 2+\log M_1
\end{array}
\]
and $S_p$ is the one in (\ref{S_p}). 
%Therefore the same argument as {\bf Case1} of the proof of Theorem \ref{thm_upper1} is valid.
%The difference appears only in finding $(x_0,t_0)\in D_2$ with $K_2(x_0,t_0)>0$.
%Let
Let us fix $(x_0,t_0) \in D_\infty$ such that 
\[
t_0=2x_0\ \mbox{and} \ t_0\ge16R^2.
\]
Then, we obtain
\[
\begin{array}{ll}
\d\left(\frac{2+t_0}{2(1+2R)}\right)^{-b}\log^{p-1}\left(\frac{2+t_0}{2(1+2R)}\right) \\
\ge \d(1+2R)^b2^{b-(p-1)}t_0^{-b}\log^{p-1}(2+t_0).
\end{array}
\]
Hence $K_7(t_0/2,t_0)>0$ follows from
\[
\begin{array}{ll}
\d \psi_2(t_0) &= t_0^{-b}\log^{p-1}(2+t_0) \\
&\d >(1+2R)^{-b}2^{-b+p-1}C_7^{-1}2^{2(p-1)S_p}(C_{f,2})^{1-p}\e^{-p(p-1)},
\end{array}
\]
which is the desired estimate as before.
This completes the proof for $a=0$ and $b<0$.
\hfill$\Box$
\vskip10pt
\par\noindent
{\bf Case 3. $\v{a<0}$ and $\v{b<0}$.}

Since $J'$ in (\ref{J'}) is estimated from below by
\[
\begin{array}{lll}
&\d\int_0^R \frac{f(-\be)^p}{(1+\be)^{1+b}} d\be \int_\be^{t-x} \frac{d\al}{(1+\al)^{1+a}}\\
&\ge\d \frac{1}{1+t-x} \int_0^R\frac{f(-\be)^p}{(1+\be)^{1+b}}d\be  \int_\be^{t-x}(\al-\be)^{-a}d\al\\
&\ge\d \frac{(t-x-R)^{1-a}}{(1-a)(1+t-x)}  \int_0^R\frac{f(-\be)^p}{(1+\be)^{1+b}}d\be,
\end{array}
\]
it follows from (\ref{frame4}) and (\ref{positive_zero_add}) that
\begin{equation}
\label{frame2_3}
\begin{array}{ll}
 \d u(x,t)& \d\ge C_0
\int_R^{t-x}\frac{d\al}{(1+\al)^{1+a}}\int_R^\al \frac{|u(y,s)|^p}{(1+\be)^{1+b}}d\be \\
 &\quad\d+C_{f,3}\e^p \frac{(t-x-R)^{1-a}}{1+t-x}
 \end{array}
\end{equation}
for $(x,t)\in D$, where
\begin{equation}
\label{Cf3}
C_{f,3}:=\frac{C_0}{2^p(1-a)}\int_0^R\frac{f(-\be)^p}{(1+\be)^{1+b}}d\be>0.
\end{equation}
Assume that an estimate
\begin{equation}
\label{case_a<0b<0g=0}
u(x,t)\ge M_n\left\{\frac{(t-x-R)^{1-b}}{1+t-x}\right\}^{a_n}\left\{\frac{(t-x-R)^{1-a}}{1+t-x}\right\}^{b_n}
\quad\mbox{in}\ D
\end{equation}
holds, where $a_n \ge 0$, $b_n>0$ and $M_n>0$. 
The sequences $\{a_n\},\{b_n\}$ and $\{M_n\}$ are defined later. 
Then it follows from (\ref{frame2_3}), (\ref{case_a<0b<0g=0}) and (\ref{Cal1}) that
\[
\begin{array}{llll}
u(x,t)&\ge\d C_0 M_n^p\int_{R}^{t-x}\frac{d\al}{(1+\al)^{1+a}}\int_{R}^{\al}\frac{(\beta-R)^{p\{a_n(1-b)+b_n(1-a)\}}}{(1+\be)^{1+b+p(a_n+b_n)}}d\be\\
&\ge\d \frac{C_0M_n^p}{(1+t-x)^{2+p(a_n+b_n)}}\int_{R}^{t-x}(\al-R)^{-a}d\al  \\
& \d \quad \times\int_{R}^{\al}(\beta-R)^{p\{a_n(1-b)+b_n(1-a)\}-b}d\be\\
&\ge\d \frac{C_0M_n^p}{\{(1-a)(pb_n+1)+(1-b)(pa_n+1)\}^2} \\
&\d \quad \times\frac{(t-x-R)^{(1-b)(pa_n+1)+(1-a)(pb_n+1)}}{(1+t-x)^{pa_n+1+pb_n+1}}
\end{array}
\]
in $D$. The inequality
(\ref{case_a<0b<0g=0}) holds for all $n\in\N$ provided
\begin{equation}
\label{a_nb_n}
\left\{
\begin{array}{l}
a_{n+1}=pa_n+1,\ a_1=0,\\
b_{n+1}=pb_n+1,\ b_1=1
\end{array}
\right.
\end{equation}
and
\[
M_{n+1}\le\frac{C_0M_n^p}{\{(1-a)(pb_n+1)+(1-b)(pa_n+1)\}^2}.
\]
%It is easy to see that
%\[
%a_n=\frac{p^{n-1}-1}{p-1},\ b_n=\frac{p^n-1}{p-1}\quad (n\in\N)
%\]
%which implies
%\[
%(pa_n+1)(pb_n+1)\le(pb_n+1)^2=b_{n+1}^2\le\frac{p^{2(n+1)}}{(p-1)^2}.
%\]
Hence we define $a_n$ by the one in (\ref{a_n1}), and set $b_n = a_{n+1}$. Moreover, $\{M_n\}$ is defined by 
\[
M_{n+1}=C_8 p^{-2n} M_n^p , \quad M_1 = C_{f,3} \e^p, 
\]where
$C_8:=C_0\{(p-1)/(2-a-b)p\}^2>0$. Therefore, making use of Lemma \ref{lemM_n} with $C=C_8$ and $\eta=\mu=p$, we obtain, by (\ref{case_a<0b<0g=0}), that 
\[
u(x,t)\ge C_8^{-1/(p-1)}\left\{\frac{(t-x-R)^{2-a-b}}{(1+t-x)^2}\right\}^{-1/(p-1)}\exp\left\{K_8(x,t)p^n\right\}
\quad\mbox{in}\ D,
\]where
\[
\begin{array}{ll}
K_8(x,t):=& \d \frac{1}{p-1}\log\left\{\frac{(t-x-R)^{1-b}}{1+t-x}\left(\frac{(t-x-R)^{1-a}}{1+t-x}\right)^p\right\}\\
& \d +\frac{1}{p-1}\log C_8-2S_p\log p+\log M_1.
\end{array}
\]
%Therefore the same argument as {\bf Case1} of the proof of Theorem \ref{thm_upper1} is valid.
%The difference appears only in finding $(x_0,t_0)\in D$ with $K_3(x_0,t_0)>0$.
%Let
Let us fix $(x_0,t_0) \in D$ such that 
\[
t_0=2x_0\ \mbox{and} \ t_0\ge4R.
\]
Then we obtain
\[
\frac{(t_0/2-R)^{1-b}}{1+t_0/2}\left(\frac{(t_0/2-R)^{1-a}}{1+t_0/2}\right)^p \ge \frac{1}{4^{(1-b)+p(1-a)}}t_0^{-b-ap},
\]
Hence $K_8(t_0/2,t_0)>0$ follows from
\[
t_0^{-b-ap}>4^{(1-b)+p(1-a)}C_8^{-1}p^{2(p-1)S_p}(C_{f,3})^{1-p}\e^{-p(p-1)}
\]
which is the desired estimate as before.
This completes the proof for $a<0$ and $b<0$.
\hfill$\Box$
\vskip10pt
\par\noindent
{\bf Case 4. $\v{a<0}$ and $\v{b=0}$.}

For $t-x\ge R$, $J'$ in (\ref{J'}) is estimated from below by
\begin{equation}
\label{a<0est1}
\begin{array}{lll}
&&\d \int_0^R \frac{f(-\be)^p}{1+\be} d\be \int_\be^{t+x} \frac{d\al}{(1+\al)^{1+a}}\\
&&\ge\d \frac{1}{1+t+x} \int_0^R\frac{f(-\be)^p}{1+\be}d\be \int_\be^{t+x} (\al-\be)^{-a}d\be\\
&&\ge\d \frac{1}{(1-a)(1+t+x)} \int_0^R\frac{f(-\be)^p(t+x-\be)^{-a}}{1+\be}d\be \\
&&\ge\d \frac{(2x)^{1-a}}{(1-a)(1+t+x)}\int_0^R\frac{f(-\be)^p}{1+\be}d\be.
\end{array}
\end{equation}
It follows from (\ref{positive_zero_add}) and (\ref{frame3}) that 
\[
u(x,t)\ge C_0
\int_R^{t-x}\frac{d\be}{1+\be}\int_\al^{t+x} \frac{|u(y,s)|^p}{(1+\al)^{1+a}}d\al
+C_{f,4}\e^p \frac{(2x)^{1-a}}{1+t+x}
\]
for $(x,t)\in D$, where
\[
C_{f,4}:=\frac{C_0}{2^p(1-a)}\int_0^R\frac{f(-\be)^p}{1+\be}d\be>0.
\]
%we employ the same argument at {\bf Case1} of the proof of Theorem \ref{thm_upper1}. 
Assume that an estimate
\begin{equation}
\label{case_a<0b=0g=0}
u(x,t)\ge M_n\left\{\log\left(\frac{1+t-x}{1+R}\right)\right\}^{a_n}\left\{\frac{(2x)^{1-a}}{1+t+x}\right\}^{b_n}
\quad\mbox{in}\ D
\end{equation}
holds, where $a_n \ge 0$, $b_n>0$ and $M_n>0$. Making use of the same computations as Case 2 of the proof of Theorem \ref{thm_upper1}, 
(\ref{case_a<0b=0g=0}) holds for all $n\in\N$ provided $\{a_n\}$ and $\{b_n\}$ are defined as in (\ref{a_nb_n}) 
and $\{M_n\}$ satisfies
\begin{equation}
\label{a<0M_n}
M_{n+1}\le\frac{C_0M_n^p}{(1-a)(pa_n+1)(pb_n+1)}.
\end{equation}
%It is easy to see that
%\[
%a_n=\frac{p^{n-1}-1}{p-1},\ b_n=\frac{p^n-1}{p-1}\quad (n\in\N)
%\]
%which implies
%\[
%(pa_n+1)(pb_n+1)\le(pb_n+1)^2=b_{n+1}^2\le\frac{p^{2(n+1)}}{(p-1)^2}.
%\]
Hence we define $a_n$ by the one in (\ref{a_n1}) and set $b_n = a_{n+1}$. Moreover, $M_n$ is defined by 
\[
M_{n+1}=C_9 p^{-2n} M_n^p , \quad M_1 = C_{f,4} \e^p,
\]where $C_9:=C_0(p-1)^2/(1-a)p^2 >0$.
Therefore, making use of Lemma \ref{lemM_n} with $C=C_9$ and $\eta=\mu=p$, we obtain, by (\ref{case_a<0b=0g=0}), that 
\[
u(x,t)\ge C_9^{-1/(p-1)}\left\{\frac{(2x)^{1-a}}{1+t+x}\log\left(\frac{1+t-x}{1+R}\right)\right\}^{-1/(p-1)}\exp\left\{K_9(x,t)p^n\right\}
\]
in $D$, where
\[
\begin{array}{ll}
K_9(x,t):=& \d \frac{1}{p-1}\log\left\{\left(\frac{(2x)^{1-a}}{1+t+x}\right)^p\log\left(\frac{1+t-x}{1+R}\right)\right\}\\
& \d +\frac{1}{p-1}\log C_9-2S_p\log p+\log M_1.
\end{array}
\]
%Therefore the same argument as {\bf Case1} of the proof of Theorem \ref{thm_upper1} is valid.
%The difference appears only in finding $(x_0,t_0)\in D$ with $K_4(x_0,t_0)>0$.
%Let
Let us fix $(x_0,t_0) \in D$ such that 
\[
t_0=2x_0\ \mbox{and} \ t_0\ge16R^2.
\]
%\[
%t_0=2x_0\quad\mbox{and}\quad t_0\ge4R^2.
%\]h
Then we obtain
\[
\left(\frac{t_0^{1-a}}{1+3t_0/2}\right)^p\log\left(\frac{2+t_0}{2+2R}\right) \ge \frac{1}{2^{p+1}}t_0^{-ap}\log(2+t_0),
\]
Hence $K_9(t_0/2,t_0)>0$ follows from
\[
\psi_1(t_0) = t_0^{-ap}\log(2+t_0) > 2^{p+1}C_9^{-1}p^{2(p-1)S_p}(C_{f,4})^{1-p}\e^{-p(p-1)}.
\]
Therefore, we obtain the desired estimate as before.
This completes the proof for $a<0$ and $b=0$.
\hfill$\Box$
\vskip10pt
\par\noindent
{\bf Case 5. $\v{a<0}$ and $\v{b>0}$.}

By virtue of (\ref{a<0est1}), (\ref{positive_zero_add}) and (\ref{frame3}), we have
\[
u(x,t)\ge C_0
\int_R^{t-x}\frac{d\be}{(1+\be)^{1+b}}\int_\al^{t+x} \frac{|u(y,s)|^p}{(1+\al)^{1+a}}d\al
+C_{f,3}\e^p \frac{(2x)^{1-a}}{1+t+x}
\]
for $(x,t)\in D$, where $C_{f,3}$ is defined in (\ref{Cf3}).
%\[
%C_f:=\frac{C_0}{2^p(1-a)}\int_0^R\frac{f(-\beta)^p}{(1+\be)^{1+b}}d\beta>0.
%\]
%we employ the same argument at {\bf Case1} of the proof of Theorem \ref{thm_upper1}. 
Assume that an estimate
\begin{equation}
\label{case_a<0b>0g=0}
u(x,t)\ge M_n\left\{\frac{t-x-R}{(1+t-x)^{1+b}}\right\}^{a_n}\left\{\frac{(2x)^{1-a}}{1+t+x}\right\}^{b_n}
\quad\mbox{in}\ D
\end{equation}
holds, where $a_n \ge 0$, $b_n>0$ and $M_n>0$. 
Similarly to the Case 3 of the proof of Theorem \ref{thm_upper1}, 
(\ref{case_a<0b>0g=0}) holds for all $n\in\N$ provided (\ref{a_nb_n}) and (\ref{a<0M_n}). 
%It is easy to see that
%\[
%a_n=\frac{p^{n-1}-1}{p-1},\ b_n=\frac{p^n-1}{p-1}\quad (n\in\N)
%\]
%which implies
%\[
%(pa_n+1)(pb_n+1)\le(pb_n+1)^2=b_{n+1}^2\le\frac{p^{2(n+1)}}{(p-1)^2}.
%\]
Hence we define $a_n$ by the one in (\ref{a_n1}), set $b_n = a_{n+1}$. Moreover, $\{M_n\}$ is defined by 
\[
M_{n+1}=C_{10} p^{-2n} M_n^p , \quad M_1 = C_{f,3} \e^p, 
\]where $C_{10}:=C_0(p-1)^2/(1-a)p^2 > 0$.
Therefore, making use of Lemma \ref{lemM_n} with $C=C_{10}$ and $\eta=\mu=p$, we obtain, by (\ref{case_a<0b>0g=0}), that 
\[
u(x,t)\ge C_{10}^{-1/(p-1)}\left\{\frac{(2x)^{1-a}}{1+t+x}\frac{t-x-R}{(1+t-x)^{1+b}}\right\}^{-1/(p-1)}\!\!\!\!\!\!\!\exp\left\{K_{10}(x,t)p^n\right\}
\quad\mbox{in}\ D,
\]where
\[
\begin{array}{ll}
K_{10}(x,t):=& \d \frac{1}{p-1}\log\left\{\left(\frac{(2x)^{1-a}}{1+t+x}\right)^p\frac{t-x-R}{(1+t-x)^{1+b}}\right\}\\
& \d +\frac{1}{p-1}\log C_{10}-2S_p\log p+\log M_1.
\end{array}
\]
%Therefore the same argument as {\bf Case1} of the proof of Theorem \ref{thm_upper1} is valid.
%The difference appears only in finding $(x_0,t_0)\in D$ with $K_5(x_0,t_0)>0$.
%Let
%\[
%t_0=x_0 + R +1 \quad\mbox{and}\quad t_0\ge6R.
%\]
Let us fix $(x_0,t_0)\in D$ such that 
\[
t_0=x_0 + R +1\ \mbox{and} \ t_0\ge4R.
\]
Then we have
\[
\frac{2^{1-a}(t_0-R-1)^{1-a}}{2t_0-R} \ge \frac{t_0^{-a}}{2},
\]
Hence $K_{10}(t_0/2,t_0)>0$ follows from
\[
t_0^{-ap}> 2(2+R)^{1+b}C_{10}^{-1}p^{2(p-1)S_p}(C_{f,3})^{1-p}\e^{-p(p-1)}.
\]
Therefore, we obtain the desired estimate as before.
This completes the proof for $a<0$ and $b>0$.
\hfill$\Box$
\vskip10pt
\par\noindent
{\bf Case 6. $\v{a=0}$ and $\v{b>0}$.}

Since $J'$ in (\ref{J'}) is estimated from below by
\[
\int_0^R \frac{f(-\be)^p}{(1+\be)^{1+b}} d\be \int_{t-x}^{t+x} \frac{d\al}{1+\al}
= \log\left(\frac{1+t+x}{1+t-x}\right)\int_0^R\frac{f(-\be)^p}{(1+\be)^{1+b}}d\be,
\]
it follows from (\ref{frame3}) and (\ref{positive_zero_add}) that
\[
\begin{array}{ll}
u(x,t)& \d \ge C_0
\int_R^{t-x}\frac{d\be}{(1+\be)^{1+b}}\int_\be^{t+x} \frac{|u(y,s)|^p}{(1+\al)^{1+a}}d\al \\
& \d \quad+C_{f,2}\e^p \log\left(\frac{1+t+x}{1+t-x}\right)
\end{array}
\]
for $(x,t)\in D$, where $C_{f,2}$ is defined in (\ref{Cf2}).
%\[
%C_f:=\frac{C_0}{2^p}\int_0^R\frac{f(-\beta)^p}{(1+\be)^{1+b}}d\beta>0.
%\]
%we employ the same argument at {\bf Case1} of the proof of Theorem \ref{thm_upper1}. 
Assume that an estimate
\begin{equation}
\label{case_a=0b>0g=0}
u(x,t)\ge M_n\left\{\frac{t-x-R}{(1+t-x)^{1+b}}\right\}^{a_n}\left\{\log\left(\frac{1+t+x}{1+t-x}\right)\right\}^{b_n}
\quad\mbox{in}\ D
\end{equation}
holds, where $a_n \ge 0$, $b_n>0$ and $M_n>0$. 
Similarly to the Case 4 of the proof of Theorem \ref{thm_upper1}, 
(\ref{case_a=0b>0g=0}) holds for all $n\in\N$ provided (\ref{a_nb_n}) 
and
\[
M_{n+1}\le\frac{C_0M_n^p}{(pa_n+1)(pb_n+1)}.
\]
%It is easy to see that
%\[
%a_n=\frac{p^{n-1}-1}{p-1},\ b_n=\frac{p^n-1}{p-1}\quad (n\in\N)
%\]
%which implies
%\[
%(pa_n+1)(pb_n+1)\le(pb_n+1)^2=b_{n+1}^2\le\frac{p^{2(n+1)}}{(p-1)^2}.
%\]
Hence we define $a_n$ as the one in (\ref{a_n1}) and set $b_n = a_{n+1}$. Moreover, $\{M_n\}$ is defined by 
\[
M_{n+1}=C_{11} p^{-2n} M_n^p , \quad M_1 = C_{f,2} \e^p, 
\]where $C_{11}:=C_0(p-1)^2/p^2 >0$.
Therefore, making use of Lemma \ref{lemM_n} with $C=C_{11}$ and $\eta=\mu=p$, we obtain, by (\ref{case_a<0b>0g=0}), that 
\[
u(x,t)\ge C_{11}^{-1/(p-1)}\left\{\frac{t-x-R}{(1+t-x)^{1+b}}\log\left(\frac{1+t+x}{1+t-x}\right)\right\}^{-1/(p-1)}\!\!\!\!\!\!\exp\left\{K_{11}(x,t)p^n\right\}
\]
in $D$, where
\[
\begin{array}{ll}
K_{11}(x,t):=& \d \frac{1}{p-1}\log\left\{\left(\log\left(\frac{1+t+x}{1+t-x}\right)\right)^p\frac{t-x-R}{(1+t-x)^{1+b}}\right\}\\
& \d +\frac{1}{p-1}\log C_{11}-2S_p\log p+\log M_1.
\end{array}
\]
%Therefore the same argument as {\bf Case1} of the proof of Theorem \ref{thm_upper1} is valid.
%The difference appears only in finding $(x_0,t_0)\in D$ with $K_6(x_0,t_0)>0$.
%Let
%\[
%t_0=x_0 + R +1 \quad\mbox{and}\quad t_0\ge4(R+1)^2.
%\]
Let us fix $(x_0,t_0)\in D$ such that 
\[
t_0=x_0 + R +1\ \mbox{and} \ t_0\ge(2+R)^2.
\]
Then we obtain
\[
\log\left(\frac{2t_0-R}{2+R}\right) \ge \log\left(\frac{t_0}{2+R}\right) \ge \frac{1}{2} \log t_0.
\]
Hence $K_{11}(t_0/2,t_0)>0$ follows from
\[
\log t_0>\{ 2^p(2+R)^{1+b}C_{11}^{-1}p^{2(p-1)S_p}(C_{f,2})^{1-p}\}^{1/p}\e^{-(p-1)}.
\]
Therefore, we obtain the desired estimate as before.
This completes the proof for $a=0$ and $b>0$.
\hfill$\Box$
\vskip10pt
\par\noindent
{\bf Case 7. $\v{a=0}$ and $\v{b=0}$.}

Since $J'$ in (\ref{J'}) is estimated from below by
\[
\int_0^R \frac{f(-\be)^p}{1+\be} d\be \int_\be^{t-x} \frac{d\al}{1+\al}
\ge \log \left(\frac{1+t-x}{1+R}\right) \int_0^R\frac{f(-\be)^p}{1+\be}d\be,
\]
it follows from (\ref{frame2}) and (\ref{positive_zero_add}) that
\begin{equation}
\label{frame2_7}
u(x,t)\ge C_0
\int_R^{t-x}\frac{d\al}{1+\al}\int_R^\al \frac{|u(y,s)|^p}{1+\be}d\be
+C_{f,5}\e^p \log \left(\frac{1+t-x}{1+R}\right)
\end{equation}
for $(x,t)\in D$, where
\[
C_{f,5}:=\frac{C_0}{2^p}\int_0^R\frac{f(-\beta)^p}{1+\be}d\beta>0.
\]
Assume that an estimate
\begin{equation}
\label{case_a=b=0g=0}
u(x,t)\ge M_n\left\{\log\left(\frac{1+t-x}{1+R}\right)\right\}^{a_n}
\quad\mbox{in}\ D
\end{equation}
holds, where $a_n \ge 0$ and $M_n>0$. 
The sequences $\{a_n\}$ and $\{M_n\}$ are defined later. 
Employing the same calculations at Case 5 of the proof of Theorem \ref{thm_upper1} and (\ref{frame2_7}), 
(\ref{case_a=b=0g=0}) holds for all $n\in\N$ provided
\[
a_{n+1}=pa_n+2,\ a_1=1
\]
and
\[
M_{n+1}\le\frac{C_0M_n^p}{(pa_n+2)^2}.
\]
It is easy to see that
\[
a_n=\frac{p^n+p^{n-1}-2}{p-1},\quad (n\in\N)
\]
which implies
\[
(pa_n+2)^2=a_{n+1}^2\le\frac{(p+1)^2}{(p-1)^2}p^{2n}.
\]
Hence $M_n$ is defined by 
\[
M_{n+1}=C_{12} p^{-2n} M_n^p , \quad M_1 = C_{f,5} \e^p,
\]where $C_{12}:=C_0(p-1)^2/(p+1)^2>0$.
Therefore, we obtain, by (\ref{case_a<0b>0g=0}), that 
\[
u(x,t)\ge C_{12}^{-1/(p-1)}\left\{\log \left(\frac{1+t-x}{1+R}\right)\right\}^{-2/(p-1)}\exp\left\{K_{12}(x,t)p^n\right\}
\]
in $D$, where
\[
\begin{array}{ll}
K_{12}(x,t):=& \d \frac{1}{p-1}\log\left\{\log^{p+1}\left(\frac{1+t-x}{1+R}\right)\right\}\\
& \d +\frac{1}{p-1}\log C_{12}-2S_p\log p+\log M_1.
\end{array}
\]
%Therefore the same argument as {\bf Case1} of the proof of Theorem \ref{thm_upper1} is valid.
%The difference appears only in finding $(x_0,t_0)\in D$ with $K_7(x_0,t_0)>0$.
%Let
%\[
%t_0=2x_0 \quad\mbox{and}\quad t_0\ge4R^2.
%\]
Let us fix $(x_0,t_0) \in D$ such that 
\[
t_0=2x_0\ \mbox{and} \ t_0\ge4(1+R)^2.
\]
Hence $K_{12}(t_0/2,t_0)>0$ follows from
\[
\log t_0>\{2^{p+1}C_{12}^{-1}p^{2(p-1)S_p}(C_{f,5})^{1-p}\}^{1/(p+1)}\e^{-p(p-1)/(p+1)}.
\]
Therefore, we obtain the desired estimate as before.
This completes the proof for $a=0$ and $b=0$.
\hfill$\Box$
\vskip10pt
\par\noindent
{\bf Case 8. $\v{a>0}$ and $\v{a+b=0}$.}

For $t-x \ge R$, $J'$ in (\ref{J'}) is estimated from below by
\[
\begin{array}{ll}
& \d \int_0^R \frac{f(-\be)^pd\be}{(1+\be)^{1+b}} \int_\be^{t-x} \frac{d\al}{(1+\al)^{1+a}} \\
& \d \ge \int_0^R \frac{f(-\be)^pd\be}{(1+\be)^{1+b}} \int_{\beta}^{R} \frac{d\al}{(1+\al)^{1+a}} \\
& \d \ge \frac{1}{(1+R)^{1+a}} \int_0^R  \frac{(R-\be)f(-\be)^p}{(1+\be)^{1+b}}d\be >0.
\end{array}
\]
Hence, we can employ the same argument as Case 6 of Theorem \ref{thm_upper1} in which the constant $C_g\e$ in (\ref{frame2}) is simply replaced with $C_{f,1}\e^p$ in (\ref{Cf1}).
Therefore the proof of Theorem \ref{thm_upper2} finishes.\hfill$\Box$
\vskip10pt

%%%%%%%%%%%%%%%%%%%%%%%%%%%%%%%%%%%%%%%%%%%%%%%%%%%%%%%%%%%%%%%%%%%%%%%%%%%%
%%%%%%%%%%%%%i%%%%%%%%%%%%%%%%%%%%%%%
%%%%%%%%%%%% Acknowledgement %%%%%%%%%%%%%%%
%%%%%%%%%%%%%%%%%%%%%%%%%%%%%%%%%%%%%
\section*{Acknowledgement}
\par
The authors would like to express their sincere gratitude to Professor Hideaki Sunagawa (Osaka Metropolitan University, JAPAN).
This work was triggered by his questions at the talk by the first author in MSJ annual meeting, Autumn 2021. 
The second author is partially supported
by the Grant-in-Aid for Scientific Research (A) (No.22H00097) and (B) (No.18H01132), 
Japan Society for the Promotion of Science. 
The third author is partially supported by the Grant-in-Aid for Young Scientist (No.20K14351), 
Japan Society for the Promotion of Science.

%%%%%%%%%%%%%%%%%%%%%%%%%%%%%%%%%%%%%%
%%%%%%%%%%%% References %%%%%%%%%%%%%%%%%%%%
%%%%%%%%%%%%%%%%%%%%%%%%%%%%%%%%%%%%%%

\bibliographystyle{plain}

\end{document}